\documentclass[12pt]{article}
\usepackage{latexsym}
\usepackage[dvips]{graphics}
\usepackage{amssymb, amsmath}
\usepackage{bm}
\usepackage{fleqn}
\renewcommand{\thefootnote}{(\arabic{footnote})}

\def\muCAS{\mu_{\mbox{\tiny CAS}}}
\def\muCLS{\mu_{\mbox{\tiny CLS}}}
\def\cP{{\cal P}}
\def\cI{{\cal I}}

\def\cK{{\cal K}}
\def\cO{{\cal O}} 
 
\def\R{\mathbb{R}}
\def\N{\mathbb{N}}
\def\F{{\sf F}} 
\def\veca{\bm{a}}

\def\vece{\bm{e}}
\def\vech{\bm{h}}
\def\vecn{\bm{n}}
\def\vect{\bm{t}}
\def\vecu{\bm{u}}
\def\vecw{\bm{w}}
\def\vecx{\bm{x}}

\def\vec0{\bm{0}}
\def\vecV{\bm{V}}
\def\vecF{\bm{F}}
\def\vecxi{\bm{\xi}}
\def\veczeta{\bm{\zeta}}
\def\eps{\varepsilon}
\def\ov{\overline} 
\def\DIV{\mathop\mathrm{div}}
\def\cosec{\mathop\mathrm{cosec}}
\def\disp{\displaystyle } 
\def\<{\langle}
\def\>{\rangle}
\def\G{\Gamma}

\makeatletter 
\long\def\@makefntext#1{\parindent 1em\noindent 
\@hangfrom{\hbox to 1.8em{\hss$^{\@thefnmark}$}}#1}
\makeatother

\def\prf{\noindent{\it Proof}.\ } 
\def\qed{{\hfill\QED}}            
\newcommand{\QED}{\hbox{\rule[0pt]{3pt}{6pt}}}
\def\barint_#1{\mathchoice
  {\mathop{\vrule width 6pt height 3 pt depth -2.5pt
    \kern -8.8pt \intop}\nolimits_{#1}}%
  {\mathop{\vrule width 5pt height 3 pt depth -2.6pt 
    \kern -6.5pt \intop}\nolimits_{#1}}%
  {\mathop{\vrule width 5pt height 3 pt depth -2.6pt
    \kern -6pt \intop}\nolimits_{#1}}%
  {\mathop{\vrule width 5pt height 3 pt depth -2.6pt
    \kern -6pt \intop}\nolimits_{#1}}}
\makeatletter
 
 \@addtoreset{equation}{section}
\makeatother

\newtheorem{Th}{Theorem}[section] 
\newtheorem{Lem}[Th]{Lemma} 
\newtheorem{Prop}[Th]{Proposition}

\newtheorem{Prob}[Th]{Problem}

\begin{document}
%
\renewcommand{\thefootnote}{\fnsymbol{footnote}}
\begin{center}{\large
Second order numerical scheme for motion of polygonal curves\\ 
with constant area speed\footnote{
The authors are supported by Czech Technical University in Prague, 
Faculty of Nuclear Sciences and Physical Engineering within the 
Jind\v rich Ne\v cas Center for 
Mathematical Modeling (Project of the Czech Ministry of Education, 
Youth and Sports LC 06052).}}\\[5pt]
\renewcommand{\thefootnote}{\arabic{footnote})}
\setcounter{footnote}{0}
\bigskip
Michal Bene\v s\footnote{
Dept. of Mathematics, 
Faculty of Nuclear Sciences and Physical Engineering, 
Czech Technical University in Prague, 
Trojanova 13, 120 00 Prague, Czech Republic. 
{\it E-mail}: benes@kmlinux.fjfi.cvut.cz}, \ 
Masato Kimura\footnote{
Faculty of Mathematics, 
Kyushu University, 
6-10-1 Hakozaki, 
Fukuoka 812-8581, Japan.
{\it E-mail}: masato@math.kyushu-u.ac.jp}
\ and \ 
Shigetoshi Yazaki\footnote{
Faculty of Engineering, 
University of Miyazaki, 
1-1 Gakuen Kibanadai Nishi, 
Miyazaki 889-2192, Japan. 
{\it E-mail}: yazaki@cc.miyazaki-u.ac.jp}
\end{center}
\bigskip
\begin{center}
\begin{minipage}[h]{0.9\linewidth}{\small
{\bf Abstract.}\ 
We study polygonal analogues of several moving boundary problems
and their time discretization which preserves the
constant area speed property.
We establish various polygonal analogues of geometric formulas for moving boundaries
and make use of the geometric formulas
for our numerical scheme and its analysis of general 
constant area speed motion of polygons.
Accuracy and efficiency of our numerical scheme
are checked through numerical simulations
for several polygonal motions such as motion by curvature 
and area-preserving advected flow etc.\\


}\end{minipage}
\end{center}

\section{\normalsize Introduction}

Polygonal analogues of several moving boundary problems
and their time discretization are investigated in this paper.
The polygonal motion is restricted within an equivalent
class of polygons.
We introduce notion of polygonal curvature,
which is consistent to polygonal analogues of
geometric variational formulas.

We propose a formulation of 
general area-preserving motion of polygonal curves by using a system of ODEs. 
The moving polygon belongs to a prescribed class of polygons, 
which is similar to the admissible class in the theory of crystalline motion by
curvature. 
There are many articles about the crystalline curvature flow and 
asymptotic behavior of solutions 
\cite{Andrews2002, GigaG2000, Girao1995, IshiiS1999, IshiwataUYY2004, 
IshiwataY2003, Yazaki2002a, Yazaki2002b, Yazaki2007a, Yazaki2007b}, etc., 
which started from the pioneer works \cite{AngenentG1989} and \cite{Taylor1993}. 
Actually, if the initial curve is a convex polygon in a crystalline
admissible class, 
then our polygonal curvature flow is nothing but the crystalline curvature flow. 
However, we consider more general polygonal moving boundary problems
in wider admissible classes of polygons.

Based on the formulation of general polygonal moving boundary problems,
we propose an implicit time discretization scheme with effective
iteration scheme for the nonlinear system in each
time step.
It has second order accuracy and preserves the constant area speed property.
In a fixed admissible class of polygons,
we prove a convergence theorem of second order for our
numerical scheme.

On the other hand, it is expected that
our polygonal analogue becomes a natural approximate
solution of a smooth moving boundary problem if
the number of edges is enough large.
It is important and interesting application of
our polygonal motion,
but we do not touch on this issue in this paper.
We only mention here that the crystalline algorithm 
and their convergence theorems for 
the motion by curvature are found in 
\cite{GigaG2000, Girao1995, GiraoK1994, IshiiS1999, Taylor1993, UY2000, UY2004} etc.

The organization of this paper is as follows. 
Fundamental notation and formulas for polygons and
polygonal motions are introduced in Section~\ref{basic}.
In Section~\ref{IVP}, a general initial value 
problem of polygonal motion in an equivalent class is 
considered, and the constant area speed condition is given.
Several basic examples of polygonal motions such as the polygonal
curvature flow and the polygonal advected flow
are also presented.
In Section~\ref{NS}, an implicit scheme of Crank-Nicolson type
and an iteration scheme for the nonlinear system in each
time step are proposed.
The proposed scheme inherits the constant area speed property
and its second order convergence is proved in 
Theorem~\ref{Th:convergence_implicit_2nd_order_scheme}.
In Section~\ref{NE}, the accuracy of our numerical scheme
is checked through various numerical simulations in comparison
with the first order explicit Euler scheme.
These simulation show that the second order scheme preserves
the constant area speed property with high accuracy.

\section{\normalsize Polygons and polygonal motions}\label{basic} 

We give basic definitions and notation 
for the polygonal motion in an equivalent class of polygons. 
In particular, the polygonal curvature is introduced
as a generalization of the crystalline curvature.
We also collect their fundamental formulas and properties
in this section.

\subsection{\normalsize Polygons}\label{polygons} 

We define a set of polygons in $\R^2$: 
\[
\cP:=\{\Gamma;\ \mbox{$\Gamma$ is a polygonal Jordan curve in $\R^2$}\}. 
\]
In this paper, we assume that any two dimensional vector
$\vecx\in \R^2$ is represented by a column vector,
and we denote its transposed row vector by $\vecx^{\rm T}$.
For $\Gamma\in\cP$, the bounded interior polygonal domain
surrounded by $\Gamma$ is denoted by $\Omega$.
For simplicity, we consider the case that
$\Omega$ is simply connected,
but many of the following arguments are valid 
in other geometrical situations.

Let $\Gamma\in \cP$ be an $N$-polygon.
The $N$ vertices of $\Gamma$ are denoted by
$\vecw_j\in \R^2$
for $j=1, 2, \ldots, N$ counterclockwise,
where $\vecw_0=\vecw_N$ and $\vecw_{N+1}=\vecw_1$.
Hereafter we use the periodic numbering convention
$\F_0=\F_N$ and $\F_{N+1}=\F_1$ 
for any quantities defined on $N$-polygon. 

For $j=1, 2, \ldots, N$, 
the $j$-th edge between $\vecw_{j-1}$ and $\vecw_j$ is defined by
\[
\Gamma_j=\{(1-\theta)\vecw_{j-1}+\theta\vecw_j;\ 0<\theta<1\},
\]
and its length is denoted by 
$|\Gamma_j|:=| \vecw_j-\vecw_{j-1}|$.
The characteristic function $\chi_j\in L^\infty(\Gamma)$ 
for $\Gamma_j$ is defined as
\[
\chi_j(\vecx):=
\left\{
\begin{array}{@{}ll@{}}
1, & \vecx\in\Gamma_j \\
0, & \vecx\in\Gamma\setminus\Gamma_j
\end{array}
\right.
\quad (j=1, 2, \ldots, N).
\]

\begin{figure}[t]
\begin{center}
\scalebox{1.2}{\includegraphics{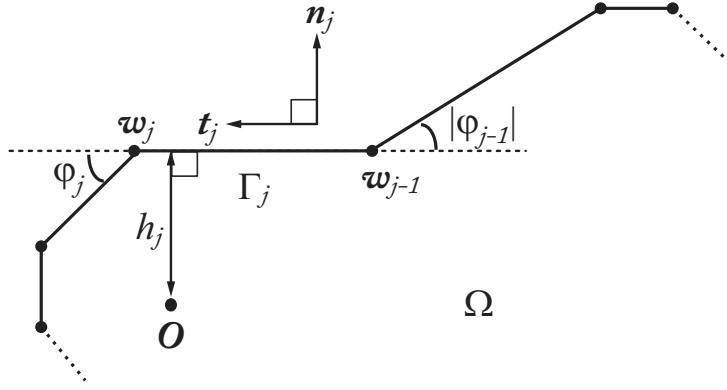}}
\caption{\small Some quantities defined on $\G_j$.}
\label{fig:polygonal-curve}
\end{center}
\end{figure}

The outward unit normal on $\Gamma_j$ is
denoted by $\vecn_j$, and the outer angle at the vertex $\vecw_j$ is denoted
by $\varphi_j\in (-\pi,\pi)\setminus \{0\}$. 
They satisfy
$\cos \varphi_j= \vecn_{j+1}\cdot\vecn_j$.
We also define the height of $\Gamma_j$ from the origin by
$h_j:=\vecw_j\cdot\vecn_j=\vecw_{j-1}\cdot\vecn_j$ 
(see Figure~\ref{fig:polygonal-curve}). 
Then the straight line 
including $\Gamma_j$ is expressed by the equation
$\vecn_j\cdot\vecx=h_j$, and
the vertices of $\Gamma\in\cP$ are given from $\{ h_j\}_j$
as
\begin{equation}\label{vecw}
\vecw_j=\left(
\begin{array}{@{}c@{}}
\vecn_j^{\rm T}\\
\vecn_{j+1}^{\rm T}
\end{array}
\right)^{-1}
\left(
\begin{array}{@{}c@{}}
h_j\\
h_{j+1}
\end{array}
\right) \quad (j=1, 2, \ldots, N).
\end{equation}
~
\begin{Prop}
Under the above conditions, we have
\begin{equation}\label{aba}
|\Gamma_j|=a_{j-1}h_{j-1}
+b_jh_j+a_jh_{j+1}
\quad (j=1, 2, \ldots, N),
\end{equation}
where 
$a_j:=\cosec\varphi_j$ and $b_j:=-\cot\varphi_{j-1}-\cot\varphi_j$. 
\end{Prop}
\prf 
We define a unit tangent vector of $\Gamma_j$ by
$\vect_j:=(\vecw_j-\vecw_{j-1})/|\Gamma_j|$.
We remark that
$\vecn_{j+1}\cdot\vect_j=-\vecn_j\cdot\vect_{j+1}=
\sin \varphi_j$
and $\vect_j\cdot\vect_{j+1}=\cos \varphi_j$.
Then, from (\ref{vecw}) and the equality:
\[
\left(
\begin{array}{@{}c@{}}
\vecn_j^{\rm T}\\
\vecn_{j+1}^{\rm T}
\end{array}
\right)
\left(
\begin{array}{@{}cc@{}}
-\vect_{j+1}
&
\vect_j
\end{array}
\right)
=
\left(
\begin{array}{@{}cc@{}}
-\vecn_j\cdot\vect_{j+1}
&0\\
0&
\vecn_{j+1}\cdot\vect_j
\end{array}
\right)
=
\frac{1}{a_j}
\left(
\begin{array}{@{}cc@{}}
1&0\\0&1
\end{array}
\right),
\]
we have
\[
\vecw_j=a_j
\left(
\begin{array}{@{}cc@{}}
-\vect_{j+1}
&
\vect_j
\end{array}
\right)
\left(
\begin{array}{@{}c@{}}
h_j\\
h_{j+1}
\end{array}
\right) 
=
-\vect_{j+1}a_jh_j+\vect_ja_jh_{j+1}.
\]
Since 
\begin{eqnarray*}
|\Gamma_j|&=&\vect_j\cdot (\vecw_j-\vecw_{j-1}) \\
&=&
\vect_j\cdot 
\left\{
(-\vect_{j+1}a_jh_j+\vect_ja_jh_{j+1})
-(-\vect_ja_{j-1}h_{j-1}+\vect_{j-1}a_{j-1}h_j)
\right\}\\
&=&
a_jh_{j+1}
-\left\{
a_j \cos \varphi_j+
a_{j-1}\cos \varphi_{j-1}
\right\}h_j
+a_{j-1}h_{j-1}\\
&=&
a_jh_{j+1}
+b_jh_j
+a_{j-1}h_{j-1},
\end{eqnarray*}
we obtain the formula (\ref{aba}).
\qed\\

The total length of $\Gamma$ is given by
\begin{equation}\label{Glength}
|\Gamma|:=\sum_{j=1}^N |\Gamma_j|
=\sum_{j=1}^N (a_j+b_j+a_{j-1})h_j
=\sum_{j=1}^N \eta_j h_j,
\end{equation}
where
$\eta_j:=a_j+b_j+a_{j-1}=\tan(\varphi_j/2)+\tan(\varphi_{j-1}/2)$. 
The area of interior domain $\Omega$ is denoted by $|\Omega|$,
which is given by 
\begin{equation}\label{area0}
|\Omega|=\frac{1}{2}\sum_{j=1}^N|\Gamma_j|h_j.
\end{equation}

The above symbols are also written as
$\vecn_j=\vecn_j(\Gamma)$,  $a_j=a_j(\Gamma)$
and $h_j=h_j(\Gamma)$ etc., 
provided we need to distinguish 
from quantities of the other polygons.

\subsection{\normalsize Equivalence classes of polygons}\label{eqc} 

For two polygons $\Gamma$ and $\Sigma\in\cP$, 
we define an equivalence relation $\Gamma\sim\Sigma$. 
We say $\Gamma\sim\Sigma$, 
if their numbers of edges are same (let it be $N$) and 
$\vecn_j(\Gamma)=\vecn_j(\Sigma)$ for all $j=1, 2, \ldots, N$ 
after choosing suitable counterclockwise numbering for $\Gamma$ and $\Sigma$. 
The equivalence class of $\Gamma\in\cP$ is denoted by 
$\cP[\Gamma]:=\{\Sigma\in\cP;\ \Sigma\sim\Gamma\}$. 

We fix an $N$-polygon $\Gamma^*\in\cP$ and let $\cP^*:=\cP[\Gamma^*]$. 
For $\Gamma$ and $\Sigma$ in $\cP^*$, 
we define the distance between them by 
\[
d(\Gamma, \Sigma):=\max_{j=1, 2, \ldots, N}|h_j(\Gamma)-h_j(\Sigma)|. 
\]
Then, it is clear that $(\cP^*, d)$ becomes a metric space, since it is 
isometrically embedded in 
$\R^N$ equipped with maximum norm $|\cdot|_\infty$ by the height function
$\vech$ defined on $\cP^*$:
\[
\vech(\Gamma):=(h_1(\Gamma), h_2(\Gamma), \ldots, h_N(\Gamma))\in\R^N 
\quad (\Gamma\in \cP^*).
\]
We assume that vectors in $\R^N$ are represented by row vectors. 
It is obvious that the image of the height function $\vech(\cP^*)$ is open in $\R^N$.

For $\Gamma\in\cP^*$ and $\eps>0$,
an $\eps$-ball in $\cP^*=\cP[\Gamma]$
with center $\Gamma$ is denoted by
\[
B(\Gamma, \eps):=
\{\Sigma\in\cP[\Gamma];\ d(\Sigma, \Gamma)<\eps\}.
\]
For an open set $\cO\subset\cP^*$ and $\Gamma\in\cO$, 
we define a positive number $\rho(\Gamma, \cO)>0$ as 
\[
\rho(\Gamma, \cO):=
\inf \{|\veca-\vech(\Gamma)|_\infty;\ 
\veca\in\R^N\setminus \vech(\cO)\}.
\]
We remark that 
$\rho(\cdot, \cO)$ is Lipschitz continuous with Lipschitz constant $1$: 
\[
|\rho(\Gamma,\cO)-\rho(\Sigma,\cO)|
\leq
d(\Gamma, \Sigma)
\quad (\Gamma, \Sigma\in\cO).
\]
For a compact set $\cK\subset\cO$, we also define
\[
\rho(\cK, \cO):=\min_{\Gamma\in\cK}\rho(\Gamma,\cO).
\]

Let $a_j=a_j (\Gamma^*)$ and $b_j=b_j(\Gamma^*)$. 
Then, from the formula (\ref{aba}), we obtain
$$
\begin{array}{@{}l}\disp
||\Gamma_j|-|\Sigma_j|| \\[5pt]\disp \qquad
=|a_{j-1}(h_{j-1}(\Gamma)-h_{j-1}(\Sigma))
+b_{j}(h_{j}(\Gamma)-h_{j}(\Sigma))
+a_{j}(h_{j+1}(\Gamma)-h_{j+1}(\Sigma))| \\[5pt]\disp\qquad
\leq C^* d(\Gamma, \Sigma)
\quad (j=1, 2, \ldots, N),
\end{array}
$$
where we define
\begin{equation}\label{C*}
C^*:=\max_{l=1, 2, \ldots, N}\{|a_{l-1}|+|b_l|+|a_l|\}.
\end{equation}

For any $\Gamma^0$ and $\Gamma^1\in\cP^*$ and for $\theta\in [0,1]$, we define 
\[
\vech^\theta:=(1-\theta)\vech (\Gamma^0)+\theta \vech (\Gamma^1)\in \R^N. 
\]
If there exists $\Gamma^\theta\in\cP^*$ with $\vech(\Gamma^\theta)=\vech^\theta$,
$\Gamma^\theta$ is called $\theta$-interpolation of $\Gamma^0$ and $\Gamma^1$.
The $\theta$-interpolation of $\Gamma^0\in\cP^*$
and $\Gamma^1\in\cP^*$ is denoted by 
$(1-\theta)\Gamma^0+\theta\Gamma^1:=\Gamma^\theta\in\cP^*$.

\subsection{\normalsize Polygonal motions}\label{polygonal} 

We consider a moving polygon $\Gamma(t)\in\cP$,
where the parameter $t$ (we call $t$ time)
belongs to an interval $\cI\subset\R$. 
For $k\in\N\cup\{0\}$, 
we call a moving polygon $\Gamma(t)$ belongs to
$C^k$-class on $\cI$, if the number of edges of $\Gamma(t)$
does not change in time and $\vecw_j\in C^k(\cI;\, \R^2)$ for all $j=1, 2, \ldots, N$.

If $k\geq 1$, we can define the normal velocity at 
$\vecx\in\Gamma_j(t)$ which is the $j$-th edge of $\Gamma(t)$. 
Let $\vecn_j(t):=\vecn_j(\Gamma (t))$.
We suppose $\vecx^*\in\Gamma_j(t^*)$ and
$\vecx^*=(1-\theta)\vecw_{j-1}(t^*)+\theta\vecw_j(t^*)$ for some $\theta\in(0,1)$, 
and define $\vecx(\theta,t):=(1-\theta)\vecw_{j-1}(t)+\theta\vecw_j(t)\in\Gamma_j(t)$. 
Then the outward normal velocity of $\Gamma_j(t^*)$ at $\vecx^*$ is defined by
\[
V_j(\vecx^*, t^*):=\dot{\vecx}(\theta, t^*)\cdot\vecn_j(t^*)
=(1-\theta)\dot{\vecw}_{j-1}(t^*)\cdot\vecn_j(t^*)
+\theta\dot{\vecw}_j(t^*)\cdot\vecn_j(t^*). 
\]
Here and hereafter, 
the (partial) derivative of $\F$ with respect to $t$ is denoted by $\dot{\F}$. 
We remark that $V_j(\cdot,t)$ is a linear function along each $\Gamma_j(t)$. 
We define the normal velocity of $\Gamma(t)$ by
\[
V(\cdot, t):=\sum_{j=1}^N V_j(\cdot, t)\chi_j(\cdot, t)\in L^\infty(\Gamma(t)),
\]
where $\chi_j(\cdot, t)\in L^\infty(\Gamma(t))$ is
the characteristic function of $\Gamma_j(t)$. 

If a moving polygon $\Gamma (t)$ belongs to a fixed equivalence class $\cP^*$ 
for all $t\in\cI$,
it is called {\it polygonal motion in} $\cP^*$ in this paper. 
Let $\vech(t)=(h_1(t),\ldots,h_N(t))\in\R^N$ be the height function
for $\Gamma (t)$.
We remark that a polygonal motion 
$\Gamma(t)$ in $\cP^*$ $(t\in\cI)$ belongs to $C^k$-class if and only if 
$\vech\in C^k(\cI;\,\R^N)$, from (\ref{vecw}).
If $\Gamma(t)$ is a $C^1$-class polygonal motion in $\cP^*$, 
its normal velocity $V_j$ of $\Gamma_j(t)$ 
is a constant on each $\Gamma_j(t)$ and it is given by $V_j(t)=\dot{h}_j(t)$. 
We denote by $\Omega(t)$ the interior domain surrounded by $\Gamma(t)$. 
\begin{Prop}
Let $\Gamma(t)$ be a $C^1$-class polygonal motion in $\cP^*$. 
Then we have
\begin{equation}\label{area2}
\frac{d}{dt}|\Omega(t)|
=
\int_{\Gamma(t)}V(\vecx,t)ds
=\sum_{j=1}^N|\Gamma_j(t)|V_j(t). 
\end{equation}
\end{Prop}
\prf 
From (\ref{aba}) and (\ref{area0}), 
we obtain
\begin{eqnarray*}
\frac{d}{dt}|\Omega(t)|
&=&
\frac{d}{dt}\left(\frac{1}{2} \sum_{j=1}^N |\Gamma_j(t)| h_j(t)\right)\\
&=& 
\frac{1}{2} \sum_{j=1}^N \left(\,
\{ a_{j-1}V_{j-1}+b_jV_j+a_jV_{j+1} \} h_j
+|\Gamma_j| V_j
\,\right)\\
&=&
\frac{1}{2} \sum_{j=1}^N V_j
\{ a_jh_{j+1}+b_jh_j+a_{j-1}h_{j-1} \}
+\frac{1}{2} \sum_{j=1}^N |\Gamma_j| V_j\\
&=&
\sum_{j=1}^N |\Gamma_j| V_j .
\end{eqnarray*}
\qed \\
%
For $\Gamma\in\cP^*$, the {\it polygonal curvature} 
$\kappa_j$ of $\Gamma_j$ is defined by
\[
\kappa_j:=\frac{\eta_j}{|\Gamma_j|},
\qquad
\eta_j:=\tan \frac{\varphi_j}{2}+\tan \frac{\varphi_{j-1}}{2}.
\]
We also define the polygonal curvature of $\Gamma$ by 
\[
\kappa:=\sum_{j=1}^N
\kappa_j\chi_j\in L^\infty(\Gamma). 
\]
The reason why this is called ``curvature'' is shown by the following proposition. 
\begin{Prop}
Let $\Gamma(t)$ $(t\in\cI)$ be a $C^1$-class polygonal motion in $\cP^*$. Then 
\[
\frac{d}{dt}|\Gamma(t)|
=\sum_{j=1}^N|\Gamma_j(t)|\kappa_j(t)V_j(t)
=\int_{\Gamma(t)}\kappa(\vecx,t)V(\vecx,t)\,ds. 
\]
\end{Prop}
\prf 
We obtain 
\[
\frac{d}{dt}|\Gamma(t)|
=\frac{d}{dt}\sum_{j=1}^{N}\eta_jh_j(t)
=\sum_{j=1}^{N}\eta_jV_j(t)
=\sum_{j=1}^N|\Gamma_j(t)|\kappa_j(t)V_j(t), 
\]
from the formula (\ref{Glength}). 
\qed

The polygonal curvature coincides with the crystalline curvature
in the crystalline motion theory (\cite{AngenentG1989, Taylor1993}). 
We, however, consider wider polygons' classes
and more general moving boundary problems.
For example, we can construct a nonconvex polygon whose all edges
have a constant positive polygonal curvature
$\kappa_1=\cdots=\kappa_N>0$ (see Section~\ref{ex4}).
We remark that such polygon is excluded in the standard crystalline theory.

\section{\normalsize Initial value problem of polygonal motion}\label{IVP} 

We consider initial value problems of polygonal motions in
an equivalent class. A general polygonal motion problem
is formulated as a system of ODEs with respect to the
height function. The notion of the constant area speed (CAS, for short)
is introduced and its necessary and sufficient condition is given.
Several concrete examples of the polygonal motion problems with
CAS property are also presented.

\subsection{\normalsize General polygonal motion problem}\label{general} 

We fix an equivalence class of $N$-polygons $\cP^*$ as in Section~\ref{polygonal}.
For an open set $\cO\subset\cP^*$ and $T_*\in(0, \infty]$, 
let $\vecF=(F_1,\ldots ,F_N)$ be a given continuous function from 
$\cO\times[0,T_*)$ to $\R^N$ with the local Lipschitz property: 
For arbitrary compact set $\cK\subset\cO$ and $T\in (0, T_*)$, 
there exists $L(\cK, T)>0$ such that 
\begin{equation}\label{Lip-cond}
|\vecF(\Gamma,t)-\vecF(\Sigma,t)|_\infty\leq L(\cK, T)\, d(\Gamma, \Sigma)
\quad (\Gamma, \Sigma\in \cK,\ t\in [0,T]). 
\end{equation}
Under the condition (\ref{Lip-cond}), 
for a compact set $\cK\subset\cO$ and $T\in (0, T_*)$, 
we also define
\[
M(\cK,T):=
\max\{|\vecF(\Gamma,t)|_\infty;\ \Gamma\in \cK,\ t\in [0,T]\}>0. 
\]
We consider the following initial value problem of polygonal motion. 
\begin{Prob}\label{PF}
For a given $N$-polygon $\Gamma^*\in\cO$, 
find a $C^1$-class polygonal motion $\Gamma(t)\in\cO$
$(0\leq t\leq T<T_*)$ such that
\[
\left\{
\begin{array}{@{}l}\disp
V_j(t)=F_j(\Gamma(t),t)
\quad (t\in [0,T],\ j=1, 2, \ldots, N)\\[5pt]\disp
\Gamma(0)=\Gamma^*.
\end{array}
\right.
\]
\end{Prob}
Under the Lipschitz condition (\ref{Lip-cond}), 
it is clear that there exists a local solution $\Gamma(t)$ in a short time interval $[0,T]$, 
since Problem~\ref{PF} can be expressed by 
an initial value problem of an ordinary differential equations 
for $\vech(t)$. 

We consider the following assumption for $F_j$:
\begin{equation}\label{areaspeed}
\sum_{j=1}^N|\Gamma_j|F_j(\Gamma,t)=\muCAS
\quad (\Gamma\in\cO, \ t\in [0,T_*)), 
\end{equation}
where $\muCAS$ is a fixed real number. 
Under the assumption (\ref{areaspeed}), 
from the formula (\ref{area2}), 
any solution $\Gamma(t)$ to Problem~\ref{PF} has the following property of 
constant area speed (CAS): 
\[
\frac{d}{dt}|\Omega(t)|=\muCAS.
\]

\subsection{\normalsize Examples of problems of polygonal motion}\label{ecmp} 

In this section, we give some examples of polygonal motions.
For several moving boundary problems for smooth curves, 
we can construct their polygonal analogues which naturally satisfy
the basic properties such as the CAS and the curve shortening (CS, for short)
properties.

\begin{Prob}[polygonal curvature flow]\label{PCF} 
For a given $N$-polygon $\Gamma^*\in\cP^*$, 
find a $C^1$-class family of $N$-polygons $\bigcup_{0\leq t\leq T}\Gamma(t)\subset\cP^*$ $(T<T_*)$ satisfying
\[
\left\{
\begin{array}{@{}l}\disp
V_j(t)=-\kappa_j(t)
\quad (t\in [0,T],\ j=1, 2, \ldots, N), \\[5pt]\disp
\Gamma(0)=\Gamma^*.
\end{array}
\right.
\]
\end{Prob}
This is a polygonal analogue of the curvature flow (curve shortening problem,
\cite{Giga2006} and see references therein). 
In the theory of crystalline motion, Problem~\ref{PCF} is 
considered in a crystalline admissible class and is called the
crystalline curvature motion.

Similar to the curvature flow for smooth curves,
the solution of Problem~\ref{PCF} has 
the CS property:
\[
\frac{d}{dt}|\Gamma(t)|
=\sum_{j=1}^N|\Gamma_j(t)|\kappa_j(t)V_j(t)
=-\sum_{j=1}^N|\Gamma_j(t)|\kappa_j(t)^2
\leq 0,
\]
and the CAS property with $\muCAS=-2\sum_{j=1}^N\tan(\varphi_j/2)$: 
\[
\frac{d}{dt}|\Omega(t)|
=-\sum_{j=1}^N|\Gamma_j(t)|\kappa_j(t)
=-\sum_{j=1}^N\eta_j
=-2\sum_{j=1}^N\tan\frac{\varphi_j}{2}
=\mbox{const.}
\]
A numerical example will be shown in Figure~\ref{fig:PCF} (left).

\begin{Prob}[area-preserving polygonal curvature flow]\label{AP-PCF} 
\quad
For a given $N$-polygon $\Gamma^*\in\cP^*$, 
find a $C^1$-class family of $N$-polygons $\bigcup_{0\leq t\leq T}\Gamma(t)\subset\cP^*$ $(T<T_*)$ satisfying
\[
\left\{
\begin{array}{@{}l}\disp
V_j(t)=\<\kappa(\cdot, t)\>-\kappa_j(t)
\quad (t\in [0,T],\ j=1, 2, \ldots, N), \\[5pt]\disp
\Gamma(0)=\Gamma^*.  
\end{array}
\right.
\]
Here $\<\kappa(\cdot, t)\>$ is the mean value of $\kappa$ on $\Gamma(t)${\rm :} 
\[
\<\kappa(\cdot, t)\>
=\frac{1}{|\Gamma(t)|}\int_{\Gamma(t)}\kappa(\vecx, t)\,ds
=\frac{\sum_{i=1}^N\eta_i}{|\Gamma(t)|}
=\frac{2\sum_{i=1}^N\tan(\varphi_i/2)}{|\Gamma(t)|}. 
\]
\end{Prob}
This is a polygonal analogue of the area-preserving curvature flow
(see \cite{Gage1986} etc.). 
Similar to the area-preserving curvature flow for smooth curves,
the solution of Problem~\ref{AP-PCF} has 
the CS property:
\[
\frac{d}{dt}|\Gamma(t)|
=\sum_{j=1}^N|\Gamma_j(t)|\kappa_j(t)V_j(t)
=-\sum_{j=1}^N|\Gamma_j(t)|(\kappa_j(t)-\<\kappa(\cdot, t)\>)^2
\leq 0,
\]
and the CAS property with $\muCAS=0$:
\[
\frac{d}{dt}|\Omega(t)|
=\<\kappa(\cdot, t)\>|\Gamma(t)|-\int_{\Gamma(t)}\kappa(\vecx, t)\,ds
=0. 
\]
Some numerical examples will be shown in Figure~\ref{fig:AP-PCF}.

In what follows, 
the mean value of $\F$ on the edge $\Gamma_j$ is denoted by 
\[
\<\F\>_j
:=\frac{1}{|\Gamma_j|}\int_{\Gamma_j}\F(\vecx)\,ds. 
\]
Let $G$ be a bounded Lipschitz domain in $\R^2$. 
We define
\[
\cO_G:=\{\Gamma\in\cP^*;\ 
\Omega(\Gamma)\supset\ov{G}\}. 
\]

\begin{Prob}[polygonal advected flow with constant area speed]\label{AP-AF} 
Let us consider $\vecu\in C^1(\R^2\setminus G;\ \R^2)$ with 
$\DIV\vecu=0$ in $\R^2\setminus \ov{G}$. 
For a given $N$-polygon $\Gamma^*\in\cO_G$, 
find a $C^1$-class family of $N$-polygons $\bigcup_{0\leq t\leq T}\Gamma(t)\subset\cO_G$ $(T<T_*)$ satisfying 
\[
\left\{
\begin{array}{@{}l}\disp
V_j(t)=\<\vecu\>_j\cdot\vecn_j
\quad (t\in [0,T],\ j=1, 2, \ldots, N), \\[5pt]\disp
\Gamma(0)=\Gamma^*.
\end{array}
\right.
\]
\end{Prob}
The solution has CAS property with $\disp\muCAS=\int_{\partial G}\vecn\cdot\vecu\,ds$\,: 
\begin{eqnarray*}
\frac{d}{dt}|\Omega(t)|
&=&
\sum_{j=1}^N|\Gamma_j(t)|\,\<\vecu\>_j\cdot\vecn_j\\
&=&
\sum_{j=1}^N \int_{\Gamma_j(t)}\vecu\cdot\vecn_j\,ds\\
&=&
\int_{\partial G}\vecu\cdot\vecn\,ds
-\int_{\Omega(t)\setminus\ov{G}}\DIV\vecu\,d\vecx\\
&=&
\int_{\partial G}\vecu\cdot\vecn\,ds,
\end{eqnarray*}
where 
$\vecn$ is the unit normal vector on $\partial G$ pointing to interior of $G$.
A numerical example will be shown in Figure~\ref{fig:AP-PAF}.

\section{\normalsize Numerical schemes}\label{NS} 

In this section, we propose an implicit time discretization of
Crank-Nicolson type
to solve the general initial value problem of polygonal motions
(Problem~\ref{PF}) and show that it preserves the CAS property 
and has a second order accuracy.
We also propose an effective iteration scheme to solve
a nonlinear system which appears in each time step.
For comparisons, we also consider an explicit Euler scheme.
We additionally give comments on the curve shortening 
and constant length speed properties
and their numerical preservation.

\subsection{\normalsize Notation}\label{NNS} 

In Section~\ref{NS}, we consider time discretization of Problem~\ref{PF} with the following notation. 
The discrete time steps are denoted by 
$0=t_0<t_1<t_2<\cdots<t_{\bar{m}}\leq T$. 
The step size which may be nonuniform and their maximum size are defined by 
\[
\tau_m:=t_{m+1}-t_m
\quad (m=0,1,\ldots,\bar{m}-1), 
\quad \tau:=\max_{0\leq m <\bar{m}}\tau_m. 
\]
Approximate solution of $\Gamma(t_m)$ is denoted by $\Gamma^m\in\cP^*$. 
Quantities of the polygon $\Gamma^m$ are denoted by 
$\vech^m=(h_1^m,\ldots,h_N^m):=(h_1(\Gamma^m),\ldots,h_N(\Gamma^m)$, 
and $\kappa_j^m:=\kappa_j(\Gamma^m)$, etc. 
We define $\vece^m:=\vech(t_m)-\vech^m\in\R^N$. 
Then we have $d(\Gamma(t_m), \Gamma^m)=|\vece^m|_\infty$. 

The discrete normal velocity $\vecV^m=(V_1^m,\ldots,V_N^m)$, 
which is an approximation of $\vecV(t_m)=\dot{\vech}(t_m)$, is defined by 
\begin{equation}\label{Vm}
\vecV^m:=\frac{\vech^{m+1}-\vech^m}{\tau_m}
\quad (m=0,1,\ldots,\bar{m}-1).
\end{equation}
Corresponding to the formula (\ref{area2}), the following formula holds. 
\begin{equation}\label{daspeed}
\frac{|\Omega^{m+1}|-|\Omega^m|}{\tau_m}
=\sum_{j=1}^N\frac{|\Gamma_j^m| +|\Gamma_j^{m+1}|}{2}V_j^m. 
\end{equation}
This has a form of sum of areas of $N$ trapezoids and is derived from (\ref{area0}) as follows: 
\begin{eqnarray*}
|\Omega^{m+1}|-|\Omega^m|
&=&
\frac{1}{2}\sum_{j=1}^N
\left( 
|\Gamma_j^{m+1}|h_j^{m+1}-|\Gamma_j^m|h_j^m
\right)\\
&=&
\frac{1}{2}\sum_{j=1}^N
\left\{
(|\Gamma_j^{m+1}|+|\Gamma_j^m|)
(h_j^{m+1}-h_j^m)
+|\Gamma_j^{m+1}|h_j^m
-|\Gamma_j^m|h_j^{m+1}
\right\}\\
&=&
\frac{\tau_m}{2}\sum_{j=1}^N
(|\Gamma_j^{m+1}|+|\Gamma_j^m|)
V_j^m
+
\frac{1}{2}\sum_{j=1}^N
\left(
|\Gamma_j^{m+1}|h_j^m
-|\Gamma_j^m|h_j^{m+1}
\right) ,
\end{eqnarray*}
where the last sum is equal to zero due to the equality (\ref{aba}).

In the following sections, 
we suppose that there exists a unique solution $\Gamma(t)$ for $0\leq t\leq T<T_*$ 
to Problem~\ref{PF} under the condition (\ref{Lip-cond}), 
and that discrete time steps $0=t_0<t_1<t_2<\cdots<t_{\bar{m}}\leq T$ 
are given a priori such as the uniform time stepping $t_m=m\tau$. 
We adopt the uniform time 
increment in the numerical examples in Section~\ref{NE}. 
It is, however, possible to apply any a posteriori adaptive time step control scheme.
Similar to the finite time extinction of the curvature flow of smooth curves, 
even in the polygonal motions, the solution polygon often has singularities 
in finite time.
For instance, 
$|\Gamma_j(t)|$ tends to zero, 
in other words, $|\kappa_j(t)|$ tends to infinity.
A posteriori adaptive time step control will be required near the blow-up time
for accurate computation.

%

\subsection{\normalsize Second order implicit scheme}\label{2IS} 

We consider the following implicit scheme for Problem~\ref{PF}.
\begin{Prob}\label{2is}
For a given $N$-polygon $\Gamma_*\in\cO$ 
and given time steps $0=t_0<t_1<t_2<\cdots<t_{\bar{m}}\leq T$, 
find polygons $\Gamma^m\in\cO$
$(m=1, 2, \ldots, \bar{m})$ such that
\[
\left\{
\begin{array}{@{}l}\disp
V_j^m=F_j(\Gamma^{m+1/2}, t_{m+1/2})
\quad (m=0, 1, 2, \ldots, \bar{m}-1,\ j=1, 2, \ldots, N), \\[5pt]
\Gamma^0=\Gamma_*, 
\end{array}
\right.
\]
where $\Gamma^{m+1/2}$ and $t_{m+1/2}$ are the $1/2$-interpolations:
\[
\Gamma^{m+1/2}:=\frac{\Gamma^m+\Gamma^{m+1}}{2}\in\cP^*, 
\quad 
t_{m+1/2}:=\frac{t_m+t_{m+1}}{2}=t_m+\frac{\tau_m}{2}.
\]
\end{Prob}
This is a generalized version of the scheme presented in
\cite{UY2004} for area-preserving crystalline curvature flow. 
\begin{Th}
We suppose the CAS property {\rm (\ref{areaspeed})}. 
Let $\Gamma^m\in\cO$ $(m=1, 2, \ldots, \bar{m})$ be a solution of Problem~\ref{2is}. 
Then it satisfies 
\[
|\Omega^{m+1}|=|\Omega^m|+\muCAS \tau_m
\quad (m=0, 1, \ldots, \bar{m}-1). 
\]
In other words, $|\Omega^m|=|\Omega(t_m)|$ holds 
if the exact solution $\Omega(t)$ of Problem~\ref{PF} exists.
\end{Th}
\prf 
Since 
$|\Gamma^{m+1/2}_j|=(|\Gamma^m_j|+|\Gamma^{m+1}_j|)/2$, 
we have
\[
\frac{|\Omega^{m+1}|-|\Omega^m|}{\tau_m}
=\sum_{j=1}^N|\Gamma^{m+1/2}_j|F_j(\Gamma^{m+1/2}, t_{m+1/2})
=\muCAS, 
\]
from the formula (\ref{daspeed}).
\qed

We remark that the numerical scheme Problem~\ref{2is}
inherits the CAS property but does not depends on the
area speed $\muCAS$. 

Since Problem~\ref{2is} is an implicit scheme, 
it is not clear whether $\Gamma^{m+1}\in\cO$ can be determined uniquely from the previous polygon 
$\Gamma^m\in\cO$, the time $t_m$, and the time step size $\tau_m$. 
Another question is how to solve the equations
\begin{equation}\label{hjm}
\vech^{m+1}=
\vech^m+\tau_m
\vecF\left(\frac{\Gamma^m+\Gamma^{m+1}}{2},\, t_{m+1/2}\right), 
\end{equation}
to obtain (approximation of) $\Gamma^{m+1}$ numerically. 
The answers to these questions will be given in Theorem~\ref{th1}
and~\ref{Th:convergence_implicit_2nd_order_scheme}.

We fix $\hat{\Gamma}\in\cO$ and $\hat{t}\in [0,T)$, which 
correspond to $\Gamma^m$ and $t_{m+1/2}$, respectively. 
Let $\cK$ be a compact convex set in $\cP^*$ with $\hat{\Gamma}\in\cK\subset\cO$. 
For $\Sigma\in\cK$ and $\hat{\tau}\in(0, \rho(\hat{\Gamma}, \cO)M(\cK, T)^{-1})$, 
we can define $\tilde{\Sigma}\in\cO$ by 
\[
\vech(\tilde{\Sigma})=
\vech(\hat{\Gamma})+\hat{\tau}
\vecF\left(\frac{\hat{\Gamma}+\Sigma}{2}, \,\hat{t}\right).
\]
We define $\Lambda(\Sigma):=
\Lambda(\Sigma; \hat{\Gamma}, \hat{t}, \hat{\tau})
:=\tilde{\Sigma}$. Then $\Lambda$ becomes a mapping from $\cK$ to $\cO$. 
We have the following lemma.
\begin{Lem}\label{contraction-lemma}
Let $\eps\in(0, \rho(\hat{\Gamma}, \cO))$ 
and $\lambda\in(0,1)$ be fixed, and let
$\hat{\cK}:=\ov{B(\hat{\Gamma}, \eps)}$. 
Suppose that $\hat{\tau}$ satisfies the condition:
\[
0<\hat{\tau}\leq\min\left\{
T-\hat{t},\ 
\frac{\eps}{M(\hat{\cK}, T)},\ 
\frac{2\lambda}{L(\hat{\cK}, T)}
\right\}.
\]
Then $\Lambda$ maps $\hat{\cK}$ into $\hat{\cK}$ and satisfies
\begin{equation}\label{dLL}
d(\Lambda(\Sigma^1), \Lambda(\Sigma^2))
\leq \lambda d(\Sigma^1,\Sigma^2)
\quad 
(\Sigma^1, \Sigma^2\in\hat{\cK}). 
\end{equation}
Namely,
$\Lambda$ is a contraction mapping on $\hat{\cK}$ and there exists a unique 
fixed point of $\Lambda$ in $\hat{\cK}$. 
\end{Lem}
\prf 
Let $\Sigma\in \hat{\cK}$. Since $\hat{\cK}$ is convex,
$(\hat{\Gamma}+\Sigma)/2\in \hat{\cK}$ holds.
We also have
\[
\left| \vech(\tilde{\Sigma})-\vech(\hat{\Gamma})\right|_\infty
=
\left| 
\hat{\tau}\ \vecF\left(\frac{\hat{\Gamma}+\Sigma}{2}, \,\hat{t}\right)
\right|_\infty
\leq
\hat{\tau} M(\hat{\cK}, T)\leq \eps .
\]
This estimate shows that $\Lambda$ is a mapping from $\hat{\cK}$ into
itself.
The estimate (\ref{dLL}) is proved as follows: 
\begin{eqnarray*}
\lefteqn{d(\Lambda(\Sigma^1), \Lambda(\Sigma^2))}\\
&&
=|\vech(\Lambda(\Sigma^1))-\vech(\Lambda(\Sigma^2))|_\infty 
=\hat{\tau}\left|\vecF\left(\frac{\hat{\Gamma}+\Sigma^1}{2},\ \hat{t}\right)
-\vecF\left(\frac{\hat{\Gamma}+\Sigma^2}{2},\ \hat{t}\right)\right|_\infty \\
&&
\leq \hat{\tau}L(\hat{\cK}, T)
d\left(\frac{\hat{\Gamma}+\Sigma^1}{2},\ \frac{\hat{\Gamma}+\Sigma^2}{2}\right) \\
&&
=\hat{\tau}L(\hat{\cK}, T)
\left|\frac{\vech(\hat{\Gamma})+\vech(\Sigma^1)}{2}
-\frac{\vech(\hat{\Gamma})+\vech(\Sigma^2)}{2}\right|_\infty \\
&&
=\frac{\hat{\tau}}{2}L(\hat{\cK}, T)\,
d\left(\Sigma^1, \Sigma^2\right)
\leq
\lambda \,d\left(\Sigma^1, \Sigma^2\right).
\end{eqnarray*}
\qed \\
The following theorem gives us an efficient numerical scheme to obtain $\Gamma^{m+1}$.
The proof is clear from Lemma~\ref{contraction-lemma}.
\begin{Th}\label{th1}
Let $\cK$ be a compact set in $\cO$ and let $\eps\in(0, \rho(\cK, \cO))$. 
We define
\[
\cK_\eps:=\ov{\bigcup_{\Sigma\in\cK}B(\Sigma, \eps)}. 
\]
For fixed $m$ $(<\bar{m})$ in Problem~\ref{2is}, 
we assume that $\Gamma^m\in\cK$ and 
\[
\tau_m\leq 
\min\left\{
\frac{\eps}{M(\cK_\eps, T)}, 
\frac{2\lambda}{L(\cK_\eps, T)}
\right\},
\]
where $\lambda\in(0,1)$. 
Then there exists uniquely 
$\Gamma^{m+1}\in\ov{B(\Gamma^m, \eps)}$ satisfying {\rm (\ref{hjm})}.

Furthermore, $\Gamma^{m+1}$ is a fixed point of 
the contraction $\Lambda_m:=\Lambda(\cdot\,; \Gamma^m, t_m, \tau_m)$ 
in $\ov{B(\Gamma^m, \eps)}$, and is given by 
the limit of $\Lambda_m^\nu(\Gamma^m)$ as $\nu\to\infty$ 
with the following estimate
\[
d(\Gamma^{m+1}, \Lambda_m^\nu(\Gamma^m))
\leq \lambda^\nu d(\Gamma^{m+1}, \Gamma^m)
\quad 
(\nu\in\N).
\]
\end{Th}

From this theorem, $\Lambda_m^\nu(\Gamma^m)$ for sufficiently large $\nu$
gives a satisfactory approximation of $\Gamma^{m+1}$.
An iteration algorithm based on this idea will be given in Section~\ref{algo}.
The following theorem shows the second order convergence of
the implicit numerical scheme (Problem~\ref{2is}).
\begin{Th} \label{Th:convergence_implicit_2nd_order_scheme}
We suppose that $\{\Gamma(t)\}_{0\leq t\leq T}$ 
be a $C^{k+1}$-class solution of Problem~\ref{PF} 
for $k=0, 1$, or $2$. 
There exists $\delta^*>0$, $\tau^*>0$, $C>0$ and 
a non-decreasing function $\omega (a)>0$ with 
\begin{equation}\label{decayomega}
\omega(a)=\left\{
\begin{array}{@{}ll}\disp
o(a^k) & \mbox{if $k=0$ or $1$,}\\[5pt]\disp
O(a^2) & \mbox{if $k=2$,}
\end{array}
\right.
\quad \mbox{as}\ a\downarrow 0,
\end{equation}
such that, if $d(\Gamma^*, \Gamma^0)\leq \delta^*$ and $\tau\leq \tau^*$, 
then $\Gamma^m\in\cO$ $(m=1, 2, \ldots, \bar{m})$ 
are inductively determined by
the implicit scheme (Problem~\ref{2is}) and  
\[
\max_{0\leq m\leq \bar{m}}
d(\Gamma(t_m), \Gamma^m)\leq 
\omega(\tau)+C d(\Gamma(0), \Gamma^0), 
\]
holds.
\end{Th}
\prf 
We put 
$\hat{\rho}:=\rho\left(\{\Gamma(t);\ 0\leq t\leq T\}, \cO\right)$, 
and fix $\delta\in(0, \hat{\rho})$ and $\eps\in(0, \hat{\rho}-\delta)$. 
We define
\[
\begin{array}{@{}l}\disp
\cK:=\ov{\bigcup_{0\leq t\leq T}B(\Gamma(t), \delta)}, 
\quad 
\cK_\eps:=\ov{\bigcup_{\Sigma\in\cK}B(\Sigma, \eps)}, \\[5pt]\disp
L:=L(\cK_\eps,T),
\quad 
R(a):=e^{\frac{L}{2}}\left(
1-\frac{aL}{2}
\right)^{-1/a}
\quad 
(0<a<2/L), \\[5pt]\disp
p_m:=|\vece^m|_\infty+\frac{\omega (\tau)}{L}
\quad (m=0, 1, 2, \ldots, \bar{m}), 
\end{array}
\]
where a non-decreasing function $\omega(a)$ $(0<a<T)$, 
which satisfies (\ref{decayomega}), will be defined in (\ref{defomega}) later. 
Since $R(\cdot)$ is an increasing function, 
there exists $\delta^*>0$ and $\tau^*>0$ such that 
\[
R(\tau^*)^T\left(
\delta^*+\frac{\omega(\tau^*)}{L}
\right)
\leq \delta,
\quad 
\tau^* <\min\left(
\frac{\eps}{M(\cK_\eps,T)},\,
\frac{2}{L}
\right).
\]

For $m=0, 1, 2, \ldots, \bar{m}-1$, 
we will prove the following inductive conditions: 
\begin{equation}\label{induction}
\Gamma^m\in\cK, \quad p_m\leq R(\tau)^{t_m}p_0
\quad\Rightarrow\quad
{}^\exists\Gamma^{m+1}\in\cK,
\quad p_{m+1}\leq R(\tau)^{t_{m+1}}p_0. 
\end{equation}
The conditions $\Gamma^0\in\cK$ and $p_0\leq R(\tau)^{0}p_0$ for the case $m=0$ 
are obviously satisfied. 

Let us assume the conditions 
$\Gamma^m\in\cK$ and $p_m\leq R(\tau)^{t_m}p_0$ for a fixed $m$. 
Then, from Theorem~\ref{th1}, 
there exists $\Gamma^{m+1}$ uniquely in $\ov{B(\Gamma^m, \eps)}\subset\cK_\eps$, 
and we have 
\begin{eqnarray}
&& \hspace{-1truecm}
\vece^{m+1}-\vece^m
=\vech(t_{m+1})-\vech(t_m)-\tau_m \vecV^m
=\tau_m\left\{\vecxi^m+\left(\vecV(t_{m+1/2})-\vecV^m\right)\right\}, 
\label{ej-ej}\\[5pt]
&& \hspace{-1truecm}
\vecxi^m
:=\frac{\vech(t_{m+1})-\vech(t_m)}{\tau_m}-\dot{\vech}(t_{m+1/2}).
\nonumber
\end{eqnarray}

The last term of (\ref{ej-ej}) is estimated as follows. 
Since $\Gamma(t_{m+1/2})\in\cK\subset\cK_\eps$ 
and $\Gamma^{m+1/2}\in\ov{B(\Gamma^m, \eps)}\subset\cK_\eps$, 
we have
\begin{eqnarray}
\lefteqn{\left|\vecV(t_{m+1/2})-\vecV^m\right|_\infty} \nonumber\\[5pt]
&&
=\left|\vecF(\Gamma(t_{m+1/2}), t_{m+1/2})-\vecF(\Gamma^{m+1/2}, t_{m+1/2})\right|_\infty
\leq L d(\Gamma(t_{m+1/2}), \Gamma^{m+1/2}) \nonumber\\[5pt]
&&
=L \left|\vech(t_{m+1/2})-\frac{\vech^m+\vech^{m+1}}{2}\right|_\infty
=L \left|\frac12(\vece^m+\vece^{m+1})-\veczeta^m\right|_\infty, 
\label{Lee}
\end{eqnarray}
where 
\[
\veczeta^m:=\frac{\vech(t_m)+\vech(t_{m+1})}{2}-\vech(t_{m+1/2}). 
\]
Combining (\ref{ej-ej}) and (\ref{Lee}), we obtain
\begin{eqnarray*}
\lefteqn{|\vece^{m+1}|_\infty 
\leq|\vece^m|_\infty+\tau_m \,|\vecxi^m|_\infty
+\tau_mL\left|\frac12(\vece^m+\vece^{m+1})-\veczeta^m\right|_\infty} \\
&&
\leq |\vece^m|_\infty+\frac{\tau_mL}{2}\left(|\vece^{m+1}|_\infty+|\vece^m|_\infty\right)
+\tau_m(|\vecxi^m|_\infty+L|\veczeta^m|_\infty).
\end{eqnarray*}
By the Taylor expansion, we can obtain 
an non-decreasing function $\omega(a)$ ($0<a<T$) 
which satisfies the condition (\ref{decayomega}) and the inequality 
\begin{equation}\label{defomega}
|\vecxi^m|_\infty+L|\veczeta^m|_\infty\leq\omega(\tau). 
\end{equation}
Hence, we have
\[
\left(1-\frac{\tau_m L}{2}\right)|\vece^{m+1}|_\infty
\leq\left(1+\frac{\tau_m L}{2}\right)
|\vece^m|_\infty+\tau_m\omega(\tau), 
\]
and this inequality is equivalent to
\[
\left(1-\frac{\tau_m L}{2}\right)p_{m+1}
\leq\left(1+\frac{\tau_m L}{2}\right)p_m. 
\]
From the inequalities
\[
\left(1-\frac{\tau_m L}{2}\right)\geq
\left(1-\frac{\tau L}{2}\right)^{\tau_m/\tau}
\quad\mbox{and}\quad
\left(1+\frac{\tau_m L}{2}\right)
\leq e^{\tau_m L/2}, 
\]
we obtain
\[
p_{m+1}
\leq
\left(1-\frac{\tau_m L}{2}\right)^{-1}
\left(1+\frac{\tau_m L}{2}\right)p_m
\leq
R(\tau)^{\tau_m}(R(\tau)^{t_m}p_0)
\leq
R(\tau)^{t_{m+1}}p_0. 
\]
The condition $\Gamma^{m+1}\in \cK$ follows from this estimate as 
\[
|\vece^{m+1}|_\infty
\leq
p_{m+1}
\leq
R(\tau^*)^{t_{m+1}}\,p_0
\leq
R(\tau^*)^T
\left(\delta^*+\frac{\omega(\tau^*)}{L}\right)
\leq
\delta.
\]
Hence, we have proved (\ref{induction}), which 
leads us to the estimate: 
\[
|\vece^m|_\infty
\leq
R(\tau^*)^T
\left(|\vece^0|_\infty+\frac{\omega(\tau)}{L}\right)
-\frac{\omega(\tau)}{L}
\leq
R(\tau^*)^T |\vece^0|_\infty
+
\frac{R(\tau^*)^T-1}{L}\omega(\tau).
\]
The assertion of the theorem is obtained by putting $C:=R(\tau^*)^T$ 
and denoting the last term $L^{-1}(R(\tau^*)^T-1)\omega(\tau)$ 
again by $\omega(\tau)$. 
\qed


\subsection{\normalsize Euler scheme}\label{ES} 

For Problem~\ref{PF}, one of the simplest numerical scheme is
the following explicit Euler scheme:
\begin{Prob}\label{es}
For a given $N$-polygon $\Gamma_*\in\cO$ 
and time steps $0=t_0<t_1<t_2<\cdots<t_{\bar{m}}\leq T$, 
find polygons $\Gamma^m\in\cO$ $(m=1, 2, \ldots, \bar{m})$ such that 
\[
\left\{
\begin{array}{@{}l}\disp
V_j^m=F_j(\Gamma^m, t_m)
\quad
(m=0, 1, 2, \ldots, \bar{m}-1,\ j=1, 2, \ldots, N), \\[5pt]\disp
\Gamma^0=\Gamma_*. 
\end{array}
\right.
\]
\end{Prob}
\noindent
The explicit Euler scheme is simple but it has only first order accuracy. 
In particular, for polygonal motions with CAS property,
we are required to use a more accurate scheme such as Problem~\ref{2is}
in order to keep its CAS property numerically. 
Similarly to the case of the implicit scheme 
(Theorem~\ref{Th:convergence_implicit_2nd_order_scheme}),
the convergence theorem of the Euler scheme is stated as follows.
\begin{Th}
We suppose the condition {\rm (\ref{Lip-cond})} 
and that $\{\Gamma(t)\}_{0\leq t\leq T}$ 
be a $C^{k+1}$-class solution of Problem~\ref{PF} 
for $k=0$ or $1$. 
There exists $\delta^*>0$, $\tau^*>0$, $C>0$ and 
a non-decreasing function $\omega (a)>0$ with 
\begin{equation}\label{decayomega-es}
\omega(a)=\left\{
\begin{array}{@{}ll}\disp
o(1) & \mbox{if $k=0$,}\\[5pt]\disp
O(a) & \mbox{if $k=1$,}
\end{array}
\right.
\quad \mbox{as}\ a\downarrow 0,
\end{equation}
such that, if $d(\Gamma^*, \Gamma^0)\leq \delta^*$ and $\tau\leq \tau^*$, 
then $\Gamma^m\in\cO$ $(m=1, 2, \ldots, \bar{m})$ is determined by
the Euler scheme (Problem~\ref{es}) and  
satisfies the estimate
\[
\max_{0\leq m\leq \bar{m}}
d(\Gamma(t_m), \Gamma^m)\leq 
\omega(\tau)+C d(\Gamma(0), \Gamma^0).
\]
\end{Th}
\prf
We define
\[
\vecxi^m
:=\frac{\vech(t_{m+1})-\vech(t_m)}{\tau_m}-\dot{\vech}(t_m)
\quad
(m=0,\ldots,\bar{m}-1).
\]
Then, by the Taylor expansion, we are able to find 
an non-decreasing function $\omega(a)$ ($0<a<T$) 
which satisfies the condition (\ref{decayomega-es}) and the inequality 
\begin{equation}\label{defomega-es}
|\vecxi^m|_\infty\leq\omega(\tau). 
\end{equation}
We put 
$\hat{\rho}:=\rho\left(\{\Gamma(t);\ 0\leq t\leq T\}, \cO\right)$, 
and fix $\delta\in(0, \hat{\rho})$ and $\eps\in(0, \hat{\rho}-\delta)$. 
We define
\[
\begin{array}{@{}l}\disp
\cK:=\ov{\bigcup_{0\leq t\leq T}B(\Gamma(t), \delta)}, 
\quad 
\cK_\eps:=\ov{\bigcup_{\Sigma\in\cK}B(\Sigma, \eps)}, \\[5pt]\disp
L:=L(\cK,T),
\quad 
p_m:=|\vece^m|_\infty+\frac{\omega (\tau)}{L}
\quad (m=0, 1, 2, \ldots, \bar{m}). 
\end{array}
\]
There exists $\delta^*>0$ and $\tau^*>0$ such that 
\[
e^{TL}\left(
\delta^*+\frac{\omega(\tau^*)}{L}
\right)
\leq \delta,
\quad 
\tau^* \leq
\frac{\eps}{M(\cK,T)}.
\]
For $m=0, 1, 2, \ldots, \bar{m}-1$, 
we will prove the following inductive conditions: 
\begin{equation}\label{induction-es}
\Gamma^m\in\cK, \quad p_m\leq e^{t_mL}p_0
\quad\Rightarrow\quad
{}^\exists\Gamma^{m+1}\in\cK,
\quad p_{m+1}\leq e^{t_{m+1}L}p_0. 
\end{equation}
The condition for $m=0$ 
are obviously satisfied. 

Let us assume the conditions 
$\Gamma^m\in\cK$ and $p_m\leq e^{t_mL}p_0$ for a fixed $m$. 
Then, from the condition $\tau^* \leq \eps /M(\cK,T)$,
$\Gamma^{m+1}$ belongs to $\ov{B(\Gamma^m, \eps)}\subset\cK_\eps$, 
and we have 
\begin{equation}\label{ej-ej-es}
\vece^{m+1}-\vece^m
=\vech(t_{m+1})-\vech(t_m)-\tau_m \vecV^m
=\tau_m\left\{\vecxi^m+\left(\vecV(t_m)-\vecV^m\right)\right\}.
\end{equation}
Since $\Gamma(t_m)$ and $\Gamma^m$ both belong to $\cK$,
we have
\begin{equation}\label{Lee-es}
\left|\vecV(t_m)-\vecV^m\right|_\infty
=
\left|\vecF(\Gamma(t_m), t_m)-\vecF(\Gamma^m, t_m)\right|_\infty
\leq L d(\Gamma(t_m), \Gamma^m) 
=L \left|\vece^m\right|_\infty.
\end{equation}
Combining (\ref{defomega-es}), (\ref{ej-ej-es}) and (\ref{Lee-es}), 
we obtain
\[
|\vece^{m+1}|_\infty
\leq\left(1+\tau_m L \right)
|\vece^m|_\infty+\tau_m\omega(\tau), 
\]
and 
\[
p_{m+1}
\leq
\left(1+\tau_m L\right)p_m
\leq
e^{\tau_mL}(e^{t_mL}\,p_0)
=
e^{t_{m+1}L}\,p_0. 
\]
Since
\[
|\vece^{m+1}|_\infty
\leq
p_{m+1}
\leq
e^{t_{m+1}L}\,p_0
\leq
e^{TL}
\left(\delta^*+\frac{\omega(\tau^*)}{L}\right)
\leq
\delta,
\]
the condition $\Gamma^{m+1}\in \cK$ follows. 
Hence, we have proved (\ref{induction-es}), which 
leads us to the estimate: 
\[
|\vece^m|_\infty
\leq
e^{LT}
\left(|\vece^0|_\infty+\frac{\omega(\tau)}{L}\right)
-\frac{\omega(\tau)}{L}
\leq
e^{LT} |\vece^0|_\infty
+
\frac{e^{LT}-1}{L}\omega(\tau).
\]
The assertion of the theorem is obtained by putting $C:=e^{LT}$ 
and denoting the last term $L^{-1}(e^{LT}-1)\omega(\tau)$ 
again by $\omega(\tau)$. 
\qed

\subsection{\normalsize Curve shortening and constant length speed property}\label{CLSP} 
As seen in Section~\ref{ecmp}, many
moving boundary problems hold the CS property:
\[
\frac{d}{dt}|\Gamma (t)|\leq 0.
\]
From (\ref{Glength}), 
a necessary and sufficient condition for CS property is
\begin{equation}\label{CScondition}
\sum_{j=1}^N\eta_j F_j(\Gamma,t)\leq 0
\quad (\Gamma\in\cO, \ t\in [0,T_*)).
\end{equation}
Similarly to the CAS property (\ref{areaspeed}),
we can also consider constant length speed (CLS, for short) property:
\[
\frac{d}{dt}|\Gamma(t)|=\muCLS. 
\]
A necessary and sufficient condition for CLS property is
\begin{equation}\label{CLScondition}
\sum_{j=1}^N\eta_j F_j(\Gamma,t)=\muCLS
\quad (\Gamma\in\cO, \ t\in [0,T_*)).
\end{equation}
An example with CLS property is the constant speed motion:
\[
F_j(\Gamma,t)=1
\quad (j=1,\ldots,N),
\qquad \muCLS=2\sum_{j=1}^N \tan \frac{\varphi_j}{2}.
\]
Another example is the length-preserving polygonal curvature flow:
\begin{equation}\label{eq:LP-PCF}
F_j(\Gamma,t)=\frac{\sum_{i=1}^N|\Gamma_j|\,\kappa_j(\Gamma)^2}
{2\sum_{i=1}^N\tan(\varphi_i/2)}-\kappa_j(\Gamma)
\quad (j=1,\ldots,N),
\qquad \muCLS=0.
\end{equation}

It is easy to check that both 
the second order implicit scheme (Problem~\ref{2is})
and the explicit Euler scheme (Problem~\ref{ES})
inherit the CS and CLS properties.
Namely, under the condition (\ref{CScondition}), we have
\[
|\Gamma^{m+1}|\leq |\Gamma^m|
\quad (m=0, 1, \ldots, \bar{m}-1),
\]
and, under the condition (\ref{CLScondition}), we have
\[
|\Gamma^{m+1}|=|\Gamma^m|+\muCLS\tau_m
\quad (m=0, 1, \ldots, \bar{m}-1).
\]
A numerical simulation for the length-preserving polygonal curvature flow
will be shown in Figure~\ref{fig:LP-PCF}.

\section{\normalsize Numerical computation}\label{NE} 

We describe an algorithm of our second order implicit scheme 
and show some numerical results. 
In this section, $\Gamma^m\in\cP^*$ $(m=0,1,\ldots,\bar{m})$
denotes the numerical solution computed by the algorithm 
described in Section~\ref{algo}. 
All computations are performed in double precision. 

\subsection{\normalsize Algorithm}\label{algo}

We describe a numerical procedure of Problem~\ref{2is}.
We suppose that an initial $N$-polygon $\Gamma^0=\bigcup_{j=1}^N\ov{\Gamma_j^0}$
is given in a prescribed equivalence class $\cP^*$, i.e.,
$\cP^*=\cP [\Gamma^0]$.
In other words, the set of normal vectors $\{\vecn_j\}_{j=1}^N$ for $\cP^*$
and the set of heights of $\Gamma_j^0$ $\vech^0=(h_1^0,\ldots,h_N^0)\in\R^N$ 
are given.
The outer angles $\{\varphi_j\}_{j=1}^N$ and 
the quantities $\{a_j\}_{j=1}^N$, $\{b_j\}_{j=1}^N$
are computed from from $\{\vecn_j\}_{j=1}^N$.
We fix the maximum computation time $T_*$ and 
the uniform time step $\tau=T_*/\bar{m}$ with the maximum time step $\bar{m}$. 

Then $\Gamma^{m+1}\in \cP^*$ is determined successively from 
$\Gamma^{m}\in\cP^*$ at the $m$-th discrete time $t_m=m\tau$ 
for $m=0, 1, \ldots, \bar{m}-1$ as follows. 
We suppose the set of heights of $\Gamma_j^m$ $\vech^m=(h_1^m,\ldots,h_N^m)\in\R^N$ 
are given.
We can 
calculate the $j$-th vertex $\vecw_j^m$ of $\Gamma^m$ by (\ref{vecw})
$(j=1, 2, \ldots, N)$.
Our algorithm including the iteration scheme to obtain 
an approximation of $\Gamma^{m+1}$ is as follows.
\[
\begin{array}{@{}l}
(1)\  \mbox{Put $\bar{\vech}:=\vech^m$.} \\[5pt]
(2)\  \mbox{Define $\hat{\Gamma}\in\cP^*$ with
$\vech(\hat{\Gamma})=\bar{\vech}$ and put $\hat{\vech}:=\bar{\vech}$.} \\[5pt]
(3)\  \mbox{Compute $\bar{\vech}:=\vech^{m}
+\vecF(\hat{\Gamma}, t_{m+1/2})\,\tau/2$.} \\[5pt]
(4)\  \mbox{If $|\bar{\vech}-\hat{\vech}|_\infty \leq\eps/2$,
then go to step (6).} \\[5pt]
(5)\  \mbox{Go to step (2).} \\[5pt]
(6)\  \mbox{Put $\vech^{m+1}:=2\bar{\vech}-\vech^m$.} 
\end{array}
\]
We note that $\hat{\Gamma}$ and $\bar{\Gamma}$ (with $\vech(\bar{\Gamma})=\bar{\vech}$)
in step (3) correspond to
$(\Lambda_m^\nu(\Gamma^m)+\Gamma^m)/2$ 
and $(\Lambda_m^{\nu+1}(\Gamma^m)+\Gamma^m)/2$,
respectively.
In the stopping condition (4),
we choose a small parameter $\eps>0$.
In the following numerical computations,
we took $\eps=10^{-15}$. 

\subsection{\normalsize Numerical examples}

In the following examples, 
several numerical computations of the evolution of $N$-sided
polygons will be shown. 
The numerical solutions were computed until the time $T_*$ with the uniform 
time increment $\tau=T_*/\bar{m}$, where $\bar{m}$ is the maximum time step. 
The figures are depicted every $M$-th time step. 
The problems except Example 7 have the CAS property with $\muCAS$, and 
the numerical solution keeps this property with the error 
$\Delta=\max_{0\leq m<\bar{m}}|\muCAS-\muCAS^m|$, 
where 
$\muCAS^m=(|\Omega^{m+1}|-|\Omega^{m}|)/\tau$ is the $m$-th discrete area speed. 
The problem in Example 7 has CLS property with $\muCLS$, and 
the numerical solution keeps this property with the error 
$\Delta=\max_{0\leq m<\bar{m}}|\muCLS-\muCLS^m|$, 
where 
$\muCLS^m=(|\Gamma^{m+1}|-|\Gamma^{m}|)/\tau$ is the $m$-th discrete
length speed. 
The following two tables indicate the data $N$, $T_*$, $\tau$, $M$ and
$\Delta$ in each example. 
\begin{table}[bht]
\begin{center}
{\small
\begin{tabular}{c|c|c|c|c|c}\hline
& \multicolumn{3}{c}{Ex.1: Figure~\ref{fig:PCF}} 
& \multicolumn{2}{|c}{Ex.2: Figure~\ref{fig:inversePCF}} 
\\\hline
& (left) & (middle) & (right) & (left) & (right) 
\\\hline\hline
$N$ & 5 & 7 & 22 & \multicolumn{2}{c}{7} 
\\\hline
$T_*$ & 0.2801 & 0.3136 & 0.335 & 1.55 & 1.55 \\\hline 
$\tau$ & \multicolumn{3}{c|}{$10^{-6}$} & $10^{-4}$ & $10^{-7}$ 
\\\hline
$M$ & 2801 & 3136 & 3350 & 775 & 775000 \\\hline
$\Delta$ & $5.04\times10^{-10}$ & $7.62\times10^{-10}$ & $1.56\times10^{-9}$ & 
$2.87\times10^{-11}$ & $2.96\times10^{7}$ 
\\\hline
\end{tabular}
}
\end{center}
\caption{Numerical parameters and $\Delta$ for Examples 1 and 2.}
\label{table1}
\end{table}
\begin{table}[tbh]
\begin{center}
{\small
\begin{tabular}{c|c|c|c|c|c|c}\hline
& \multicolumn{2}{|c}{Ex.3: Figure~\ref{fig:AP-PCF}} 
& \multicolumn{1}{|c}{Ex.5: Figure~\ref{fig:AP-PAF}}
& \multicolumn{2}{|c}{Ex.6: Figure~\ref{fig:AP-adPCF}}
& \multicolumn{1}{|c}{Ex.7: Figure~\ref{fig:LP-PCF}}
\\\hline
&& (upper) & (lower) & (upper) & (lower) &
\\\hline\hline
$N$ & 9 & 12 & 12 & \multicolumn{2}{|c|}{32} & 18
\\\hline
$T_*$ & 7.56 & 19.4 & 20 & 10 & 3.65 & 0.27 \\\hline 
$\tau$ & \multicolumn{2}{c|}{$10^{-5}$} 
& $10^{-4}$ & \multicolumn{2}{c|}{$10^{-4}$} & $10^{-4}$ \\\hline
$M$ & 37800 & 97000 & 10000 & 5000 & 1825 & 27 \\\hline
$\Delta$ & $1.51\times10^{-9}$ & $1.07\times10^{-9}$ 
& $2.81\times10^{-7}$ & $4.26\times 10^{-9}$ & $1.42\times 10^{-10}$ & $5.11\times10^{-11}$ 
\\\hline
\end{tabular}
}
\end{center}
\caption{Numerical parameters and $\Delta$ for Examples 4--7.}
\label{table2}
\end{table}

\subsubsection{\normalsize Example 1 --- polygonal curvature flow}\label{ex1} 

\begin{figure}[ht]
\begin{center}
\scalebox{0.8}{\includegraphics{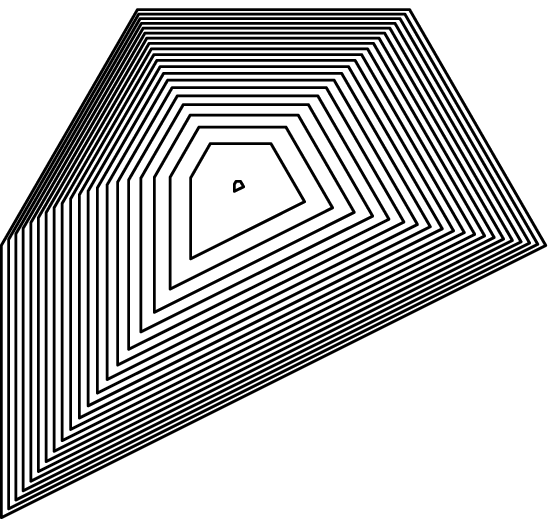}}\hspace{1truecm}
\scalebox{0.8}{\includegraphics{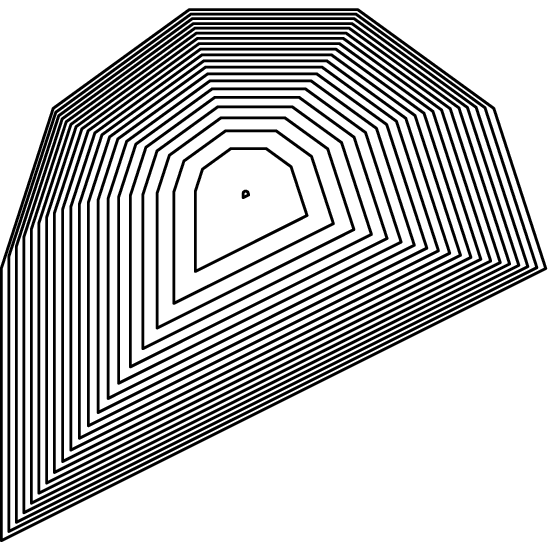}}\hspace{1truecm}
\scalebox{0.8}{\includegraphics{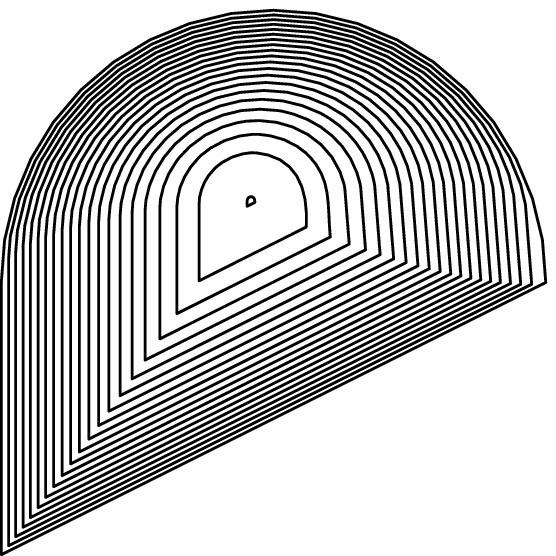}}
\caption{\small Evolution by polygonal curvature flow. 
}
\label{fig:PCF}
\end{center}
\end{figure}

Figure~\ref{fig:PCF} indicates the evolution of solution polygons to Problem~\ref{PCF}, 
starting from the initial polygon being the outermost $N$-sided polygon which is 
a combination of 
an upper half of a regular $2(N-2)$-polygon and a triangle.
Each solution polygon evolves from outside to inside
and has the CAS property with
$\muCAS=-2\sum_{j=1}^N\tan(\varphi_j/2)$. 
The numerical solutions keep the CAS property 
very accurately as shown in Table~\ref{table1}.

\subsubsection{\normalsize Example 2 --- backward polygonal curvature flow}\label{ex2} 

\begin{figure}[ht]
\begin{center}
\scalebox{0.8}{\includegraphics{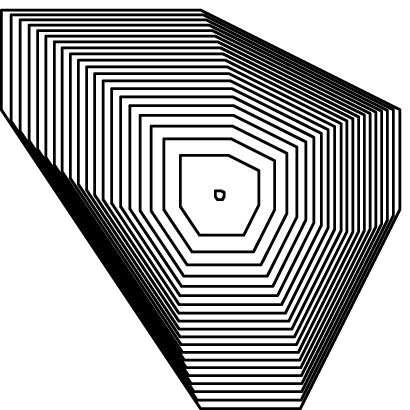}}\hspace{1truecm}
\scalebox{0.8}{\includegraphics{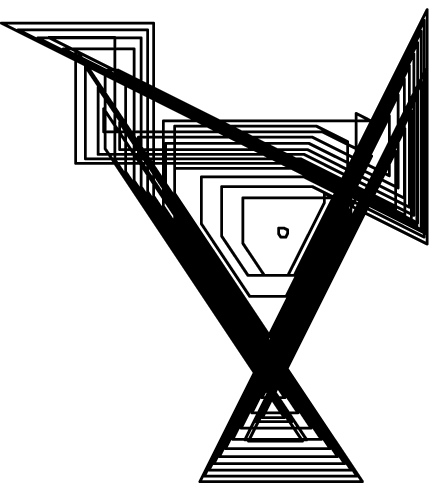}}
\caption{\small Simulations of the backward polygonal curvature flow 
by the second order scheme (left) and by the Euler method (right). 
}
\label{fig:inversePCF}
\end{center}
\end{figure}

Problem~\ref{PCF} can be computed backward in time. 
Figure~\ref{fig:inversePCF} (left) indicates the evolution of solution polygons to 
the backward polygonal curvature flow $V_j(t)=\kappa_j(t)$ ($j=1, 2, \ldots, 7$). 
The initial polygon is the innermost 7-sided polygon and 
the solution polygons evolve from inside to outside. 
The above process can be followed by our second order scheme accurately. 

We note that the backward curvature flow for smooth curves is ill-posed 
since it becomes a backward parabolic problem.
Actually, even in the case of 7-sided polygon's motion,
it is hard to compute the backward polygonal curvature flow by using 
the Euler scheme (Problem~\ref{es}).
Figure~\ref{fig:inversePCF} (right) indicates an easy breakdown of the Euler scheme 
in spite of using a smaller $\tau$ than the one in the second order scheme. 

\subsubsection{\normalsize Example 3 --- area-preserving polygonal curvature flow}\label{ex3} 

\begin{figure}[ht]
\begin{center}
\scalebox{0.6}{\includegraphics{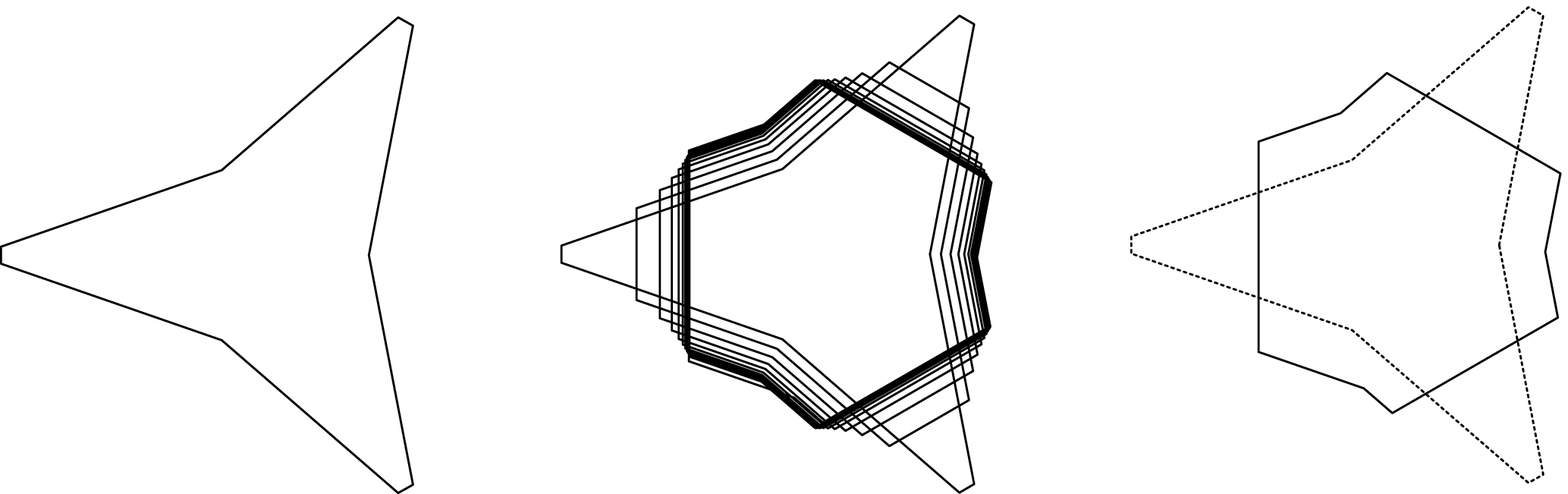}}\\[5pt]
\scalebox{0.6}{\includegraphics{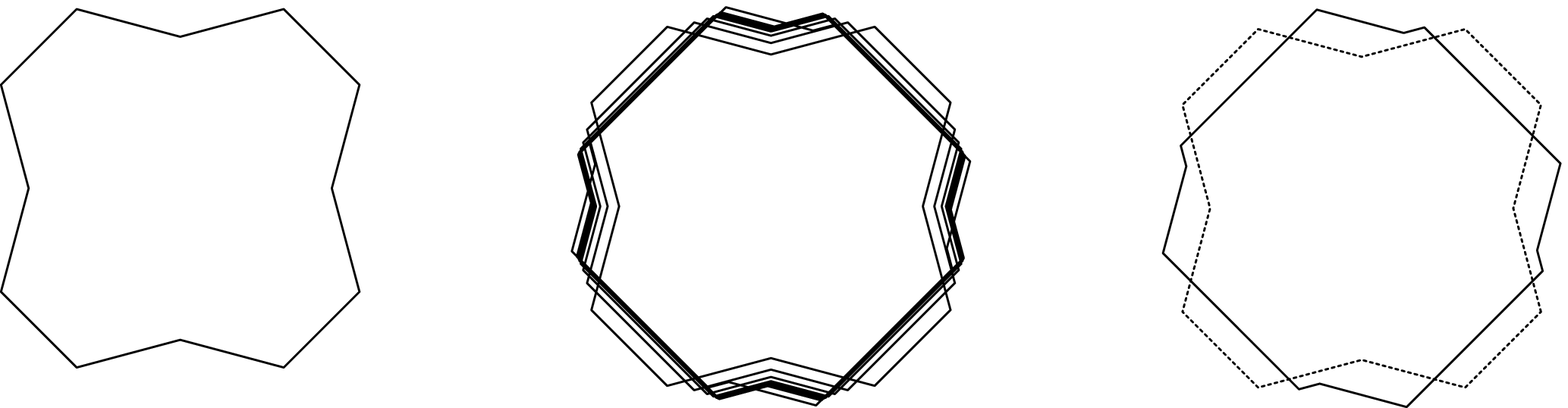}}
\caption{\small Evolutions by the area-preserving polygonal curvature flow.}
\label{fig:AP-PCF}
\end{center}
\end{figure}

Figure~\ref{fig:AP-PCF} (middle--upper/lower) 
shows two examples of polygonal motions according to Problem~\ref{AP-PCF}. 
The initial polygons are given as Figure~\ref{fig:AP-PCF} (left--upper/lower). 
Figure~\ref{fig:AP-PCF} (right--upper/lower) shows the final polygon and 
the initial polygon (dotted curve). 
The solution has CAS property with $\muCAS=0$. 

In both upper and lower examples, there exist stationary solutions
as shown in Figure~\ref{fig:AP-PCF-ss} (middle/right). 
The polygon starting from a symmetric initial shape
approaches to one of the stationary solutions and stays there for a while. 
However, since the stationary solution has a saddle-point instability,
after a while,
the polygon is drifted away from the stationary solution along the
unstable manifold and loses its symmetry.

\subsubsection{\normalsize Example 4 --- stationary solutions}\label{ex4} 

\begin{figure}[ht]
\begin{center}
\scalebox{0.54}{\includegraphics{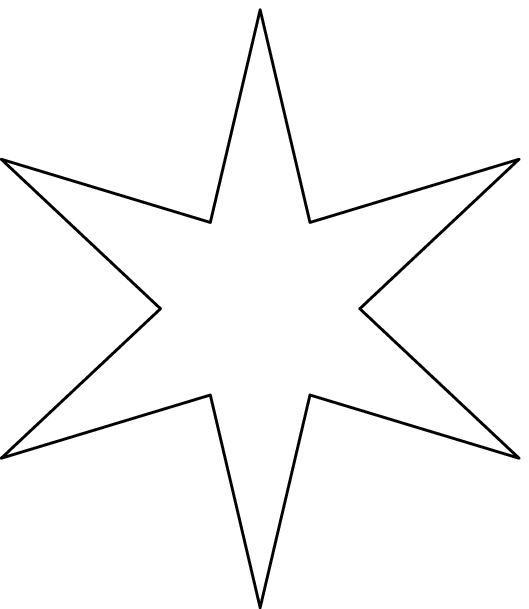}}\hspace{1truecm}
\scalebox{0.65}{\includegraphics{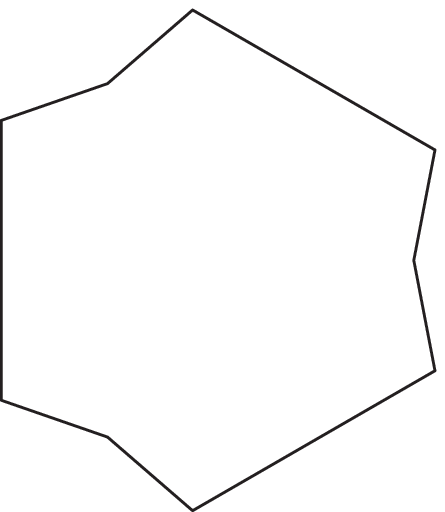}}\hspace{1truecm}
\scalebox{0.7}{\includegraphics{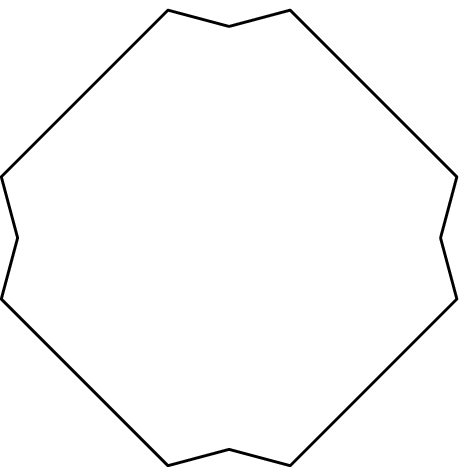}}
\caption{\small Stationary solutions of Problem~\ref{AP-PCF}.}
\label{fig:AP-PCF-ss}
\end{center}
\end{figure}

A polygon which has a constant polygonal curvature
(i.e. $\kappa_1=\cdots =\kappa_N$) is a stationary solution of 
Problem~\ref{AP-PCF}. 
Obviously, regular polygons are stationary solutions. 
Besides the regular polygons, we have infinite many stationary solutions. 
For instance, 
an $n$-fold star shaped polygon is a stationary solution 
as well as the 6-fold star (Figure~\ref{fig:AP-PCF-ss} (left)), and 
an $n$-fold non-sharp star shaped polygon is also a stationary solution 
as well as the 3-fold/4-fold non-sharp stars in Figure~\ref{fig:AP-PCF-ss} (middle/right). 
For the $n$-fold non-sharp star polygon,
there are two kinds of outer angles $\varphi_0<0$ and $\varphi_1>0$
and two kinds of edge lengths $d_1 < d_2$
with the corresponding polygonal curvatures: 
$\kappa_1=(\tan(\varphi_0/2)+\tan(\varphi_1/2))/d_1$ and 
$\kappa_2=2\tan(\varphi_1/2)/d_2$. 
We have constant polygonal curvature polygon
if $d_1/d_2=(\tan(\varphi_0/2)+\tan(\varphi_1/2))/(2\tan(\varphi_1/2))$.
The polygon in Figure~\ref{fig:AP-PCF-ss} (middle/right) belongs
the same equivalence class of Figure~\ref{fig:AP-PCF} (upper/lower). 
%
\subsubsection{\normalsize 
Example 5 --- polygonal advected flow with constant area speed}\label{ex5} 

\begin{figure}[b]
\begin{center}
\scalebox{0.6}{\includegraphics{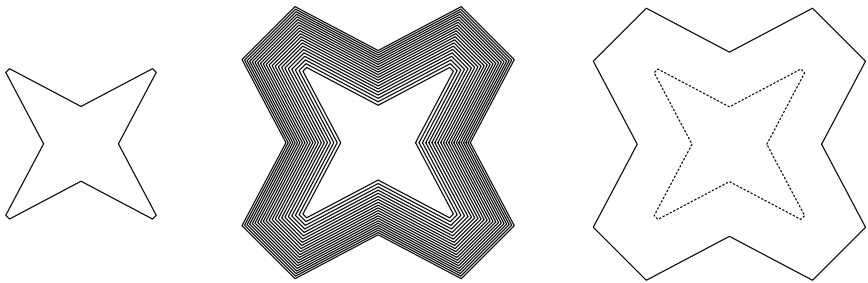}}
\caption{\small Evolution by polygonal advected flow with constant area speed.}
\label{fig:AP-PAF}
\end{center}
\end{figure}

Figure~\ref{fig:AP-PAF} (middle)
shows an example of polygonal motions according to Problem~\ref{AP-AF} 
with $\vecu (\vecx) =\vecx /(2\pi |\vecx|^2)$ which is a divergence-free
vector field defined on $\R^2\setminus\{\vec0\}$.
The initial polygons are given as Figure~\ref{fig:AP-PAF} (left)
whose center is the origin. 
Figure~\ref{fig:AP-PAF} (right) shows the final polygon and 
the initial polygon (dotted curve). 
The problem has the CAS property with $\muCAS=1$. 
The numerical solution keeps the CAS property 
very accurately as shown in Table~\ref{table2}.

\subsubsection{\normalsize 
Example 6 --- area-preserving polygonal advected-curvature flow}\label{ex6} 

\begin{figure}[ht]
\begin{center}
\scalebox{0.6}{\includegraphics{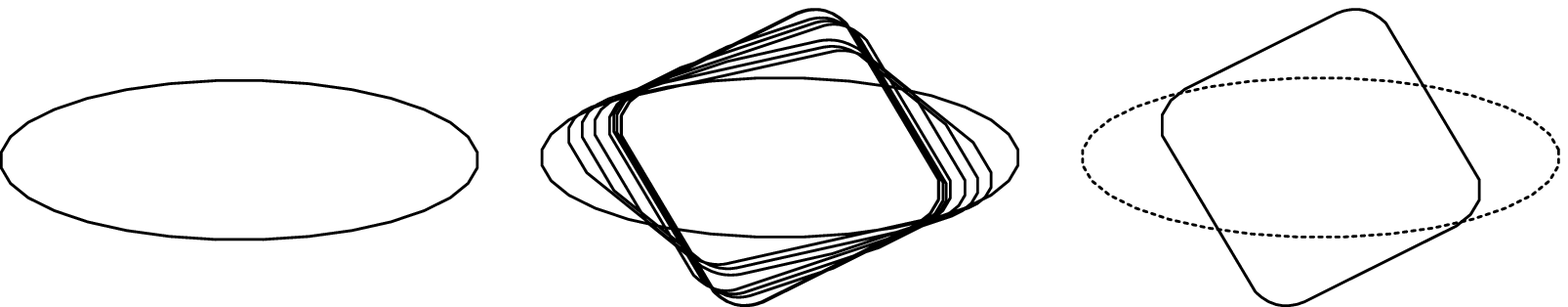}}\\[5pt]
\scalebox{0.6}{\includegraphics{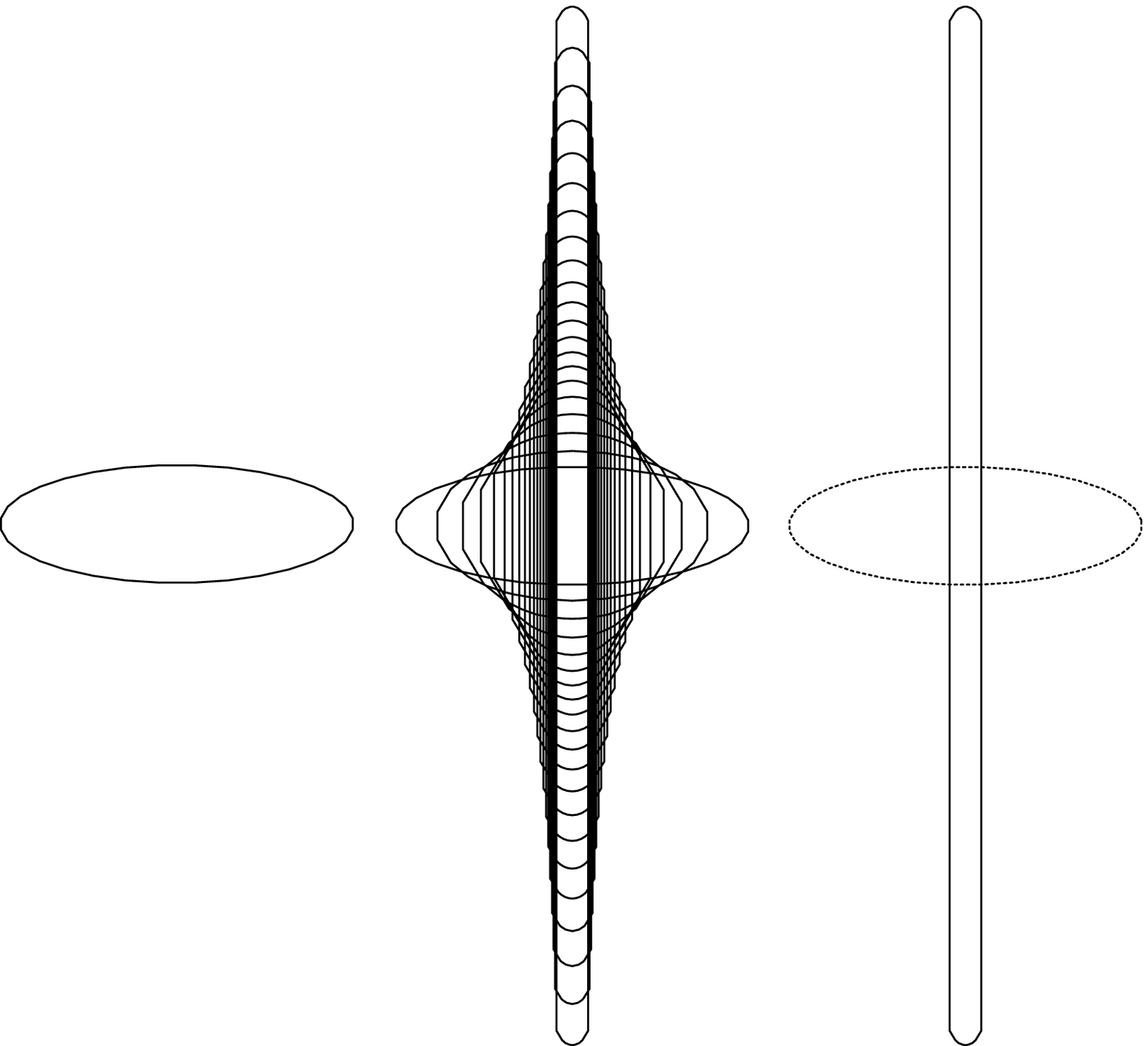}}
\caption{\small Evolution by area-preserving polygonal advected-curvature flow.}
\label{fig:AP-adPCF}
\end{center}
\end{figure}

Figure~\ref{fig:AP-adPCF} (middle--upper/lower) 
shows two examples of polygonal motions according to 
combination of Problem~\ref{AP-PCF} and Problem~\ref{AP-AF},
i.e., 
\[
V_j=\<\kappa(\cdot, t)\>-\kappa_j(t)+\<\vecu\>_j\cdot\vecn_j 
\quad 
(j=1,\ldots,N).
\]
The divergence-free vector field is given by 
$\vecu(\vecx)=x_1x_2(-x_1, x_2)$ (upper) and $\vecu(\vecx)=(-x_1, x_2)$
(lower), respectively.
The common initial polygon is given as Figure~\ref{fig:AP-adPCF}
(left--upper/lower), where the center is the origin the vertices 
are on the ellipse with ratio 3:1. 
Figure~\ref{fig:AP-adPCF} (right--upper/lower) shows the final polygon and 
the initial polygon (dotted curve). 
The both problems have CAS property with $\muCAS=0$
(area preserving), and the numerical solutions preserve their areas
very accurately as shown in Table~\ref{table2}.

\subsubsection{\normalsize 
Example 7 --- length-preserving polygonal curvature flow}\label{ex7} 

\begin{figure}[ht]
\begin{center}
\scalebox{0.7}{\includegraphics{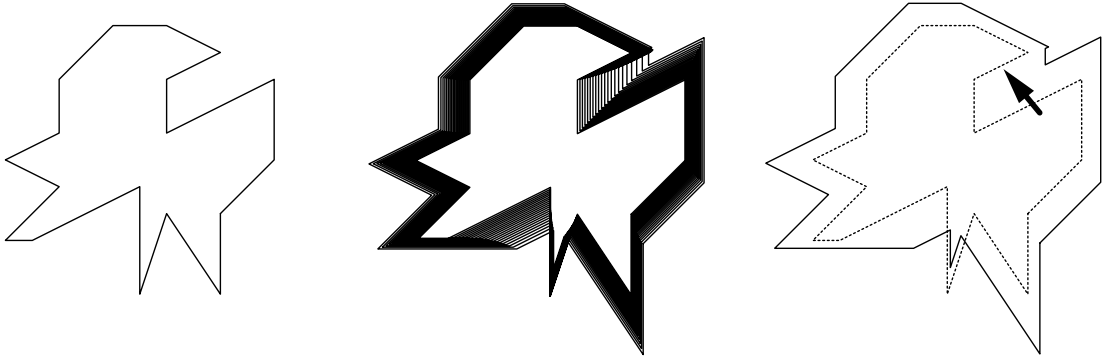}}
\caption{\small Evolution by the length-preserving polygonal curvature flow. }
\label{fig:LP-PCF}
\end{center}
\end{figure}

Figure~\ref{fig:LP-PCF} (middle) shows the evolution of solution polygons to 
length-preserving polygonal curvature flow (\ref{eq:LP-PCF}) given in
Section~\ref{CLSP}. 
The initial polygon is the 18-sided polygon in Figure~\ref{fig:LP-PCF} (left). 
Figure~\ref{fig:LP-PCF} (right) indicates the final polygon and 
the initial polygon (dotted curve). 
At the time close to $T_*$, 
the length of an edge (pointed by the arrow) tends to zero, 
and the computation stops. 
The length-preserving property ($\muCLS=0$) is numerically realized
with high accuracy  as shown in Table~\ref{table2}.


\begin{thebibliography}{XX}

\bibitem{Andrews2002}
B. Andrews, 
{\sl Singularities in crystalline curvature flows}, 
Asian J. Math. {\bf 6} (2002) 101--122. 

\bibitem{AngenentG1989}
S. Angenent and M. E. Gurtin, 
{\sl Multiphase thermomechanics with interfacial structure, 
{\rm 2}. Evolution of an isothermal interface}, 
Arch. Rational Mech. Anal. {\bf 108} (1989) 323--391.

\bibitem{Gage1986}
M. Gage, 
{\sl On an area-preserving evolution equation for plane curves}, 
D.M. DeTurck (Ed.), Nonlinear Problems in Geometry, Contemp. Math. 
{\bf 51} (1986) 51--62.

\bibitem{Giga2006}
Y. Giga, 
{\sl Surface evolution equations, 
A level set approach},
Monographs in Mathematics, 99. 
Birkhauser Verlag, Basel, (2006).

\bibitem{GigaG2000}
M.-H. Giga and Y. Giga, 
{\sl Crystalline and level set flow 
-- convergence of a crystalline algorithm for a 
general anisotropic curvature flow in the plane}, 
Free boundary problems: theory and applications, I (Chiba, 1999), 
GAKUTO Internat. Ser. Math. Sci. Appl., 
Gakk\={o}tosho, Tokyo {\bf 13} (2000) 64--79. 

\bibitem{Girao1995}
P. M. Gir\~ao,
{\sl Convergence of a crystalline algorithm for the motion 
of a simple closed convex curve by weighted curvature},
SIAM J. Numer. Anal. {\bf 32} (1995) 886--899.

\bibitem{GiraoK1994}
P. M. Gir\~ao and R. V. Kohn, 
{\sl Convergence of a crystalline algorithm for the heat
equation in one dimension and for the motion of a
graph by weighted curvature}, 
Numer. Math. {\bf 67} (1994) 41--70.

\bibitem{IshiiS1999}
K. Ishii  and H. M. Soner , 
{\sl Regularity and convergence of crystalline motion}, 
SIAM J. Math. Anal. {\bf 30} (1999) 19--37.

\bibitem{IshiwataUYY2004}
T. Ishiwata , T. K. Ushijima , H.  Yagisita and S. Yazaki, 
{\sl Two examples of nonconvex self-similar solution 
curves for a crystalline curvature flow}, 
Proc. Japan Academy {\bf 80}, Ser. A, No. 8 (2004), 151--154. 

\bibitem{IshiwataY2003}
T. Ishiwata and S. Yazaki, 
{\sl On the blow-up rate for fast blow-up solutions arising in an 
anisotropic crystalline motion}, 
J. Comp. App. Math. {\bf 159} (2003), 55--64. 

\bibitem{Taylor1993}
J. E. Taylor,
{\sl Motion of curves by crystalline curvature, 
including triple junctions and boundary points}, 
Diff. Geom.: partial diff. eqs. on manifolds 
(Los Angeles, CA, 1990), Proc. Sympos. Pure Math., 
{\bf 54} (1993), Part I, 417--438, 
AMS, Providencd, RI.

\bibitem{UY2000}
T. K. Ushijima and S. Yazaki, 
{\sl Convergence of a crystalline algorithm for the motion of
a closed convex curve by a power of curvature $V=K^\alpha$}, 
SIAM J. Numer. Anal. {\bf 37} (2000) 500--522.

\bibitem{UY2004}
T. K. Ushijima and S. Yazaki, 
{\sl Convergence of a crystalline approximation for an 
area-preserving motion}, 
Journal of Computational and Applied Mathematics {\bf 166} (2004), 427--452. 

\bibitem{Yazaki2002a}
S. Yazaki, 
{\sl Asymptotic behavior of solutions to an expanding 
motion by a negative power of crystalline curvature}, 
Adv. Math. Sci. Appl. {\bf 12} (2002), 227--243. 

\bibitem{Yazaki2002b}	
S. Yazaki, 
{\sl On an area-preserving crystalline motion}, 
Calc. Var. {\bf 14} (2002), 85--105. 

\bibitem{Yazaki2007a}
S. Yazaki, 
{\sl Motion of nonadmissible convex polygons by crystalline curvature}, 
Publications of Research Institute for Mathematical Sciences {\bf 43} (2007), 155--170. 

\bibitem{Yazaki2007b}
S. Yazaki, 
{\sl Asymptotic behavior of solutions to an area-preserving motion by crystalline curvature}, 
Kybernetika {\bf 43} (2007), 903--912.

\end{thebibliography}
\end{document}